\documentclass[10pt]{amsart}

\usepackage{graphicx}
\usepackage{xcolor}
\usepackage{amssymb}
\usepackage{bbm}
\usepackage{amsmath,mathtools,enumerate}
\usepackage{amsthm}
\usepackage{stmaryrd}
\usepackage{tikz-cd}
\usepackage{faktor}
\usepackage{hyperref}
\usepackage{adjustbox}
\usepackage{verbatim}

\newcommand{\C}{\mathbb{C}}

\newcommand{\N}{\mathbb{N}} 

\renewcommand{\O}{\mathcal{O}}
\renewcommand{\L}{\mathcal{L}}
\newcommand{\T}{\mathcal{T}}

\newcommand{\G}{\mathcal{G}}

\newcommand{\K}{\mathcal {K}}

\newcommand{\Ker}{\mathrm{Ker}} 

\newcommand{\id}{\mathrm{id}}

\newcommand{\dom}{\mathrm{dom}}
\newcommand{\ran}{\mathrm{ran}}

\newcommand{\cstar}{\texorpdfstring{C$^*$}{C*}}
\newcommand{\safeK}{\texorpdfstring{$K$}{K}}

\theoremstyle{definition}

\newtheorem{thm}{Theorem}[section]
\newtheorem{lem}[thm]{Lemma}
\newtheorem{prop}[thm]{Proposition}
\newtheorem{dfn}[thm]{Definition}
\newtheorem{cor}[thm]{Corollary}

\numberwithin{equation}{section}

\title[
Every classifiable simple C$^*$-algebra has a C$^*$-diagonal]{
Every classifiable simple C$^*$-algebra has a C$^*$-diagonal}

\author{Samuel Evington}
\address{Samuel Evington, School of Mathematics and Statistics, University College Dublin, Belfield, Dublin 4, D04 V1W8, Ireland}
\email{sam.evington@ucd.ie}

\author{Philipp Sibbel}
\address{Philipp Sibbel, Mathematical Institute, University of M\"unster, Ein\-stein\-strasse 62, 48149 M\"unster, Germany}
\email{philipp.sibbel@uni-muenster.de}

\thanks{Funded by the Deutsche Forschungsgemeinschaft (DFG, German Research Foundation) under Germany’s Excellence Strategy – EXC 2044/2 – 390685587, Mathematics Münster – Dynamics – Geometry – Structure; the Deutsche Forschungsgemeinschaft (DFG, German Research Foundation) – Project-ID 427320536 – SFB 1442; and ERC Advanced Grant 834267 - AMAREC}

\begin{document}
	
	\begin{abstract}  
        We construct C$^*$-diagonals in all UCT Kirchberg algebras. Together with known constructions of C$^*$-diagonals in the stably finite setting by Xin Li, this shows that every classifiable C$^*$-algebra (in the sense of the Elliott classification programme) has a C$^*$-diagonal. In fact, we show that every UCT Kirchberg algebra admits an untwisted principal étale groupoid model with a finite-dimensional unit space.
	\end{abstract}
	
	\maketitle
	
	\section*{Introduction}
	\numberwithin{equation}{section}
	\renewcommand{\thethm}{\Alph{thm}}

    The class of simple, separable, nuclear, non-elementary C$^*$-algebras that satisfy both the universal coefficient theorem (UCT) of Rosenberg and Schochet (\cite{RS87}) and any of the regularity conditions of the recently-established Toms--Winter conjecture (\cite{TomsWinter, Thiel26,OAI:Toms-Winter}) are classified by an invariant built out of $K$-theory and traces, thanks to the work of many hands over the last 50 years on the Elliott classification programme (see \cite{El95, Ki95, winter2018structure, WhiteICM, GL23, CGSTW} for an overview). We shall refer to these C$^*$-algebras as \emph{classifiable} in this introduction. This class of C$^*$-algebras includes the non-elementary simple approximately finite-dimensional algebras (\cite{elliott1976classification}), the irrational rotation algebras (\cite{rieffel1981RotationAlgebras}), the Cuntz algebras (\cite{cuntz1977simple}), as well as many other C$^*$-algebras constructed via inductive limits (\cite{EGLN25, GLN20a, GLN20b}) or from dynamical systems (\cite{MatuiSato12, MatuiSato14, KerrSzabo20, GGNV24}).

    The headline results of the Elliott classification programme are uniqueness theorems that, loosely speaking, state that any two classifiable C$^*$-algebras with the same invariant are isomorphic. However, it is equally important to construct examples realising all possible invariants (\cite{El96}). Such examples provide concrete descriptions of C$^*$-algebras and can reveal structure not apparent from their classification invariants. Moreover, having multiple concrete models of a given C$^*$-algebra is of great value because structure or symmetry that is immediately apparent in one model may not be in another (see \cite[Section~13]{STW25_99problems}).

    Groupoid models provide one means of realising classifiable C$^*$-algebras and make methods from topological dynamics available in the study of the resulting C$^*$-algebras. 
    By work of Xin Li in the stably finite case (\cite{li2020every}) and earlier results in the purely infinite case (\cite{spielberg2007graph,katsura2008class,yeend2006topological,li2019cartan, CFH20}), 
    every classifiable C$^*$-algebra $A$ is isomorphic to the reduced C$^*$-algebra $C_r^*(\G,\Sigma)$ of a twisted groupoid $(\G, \Sigma)$ that is topologically principal, second-countable, locally compact, Hausdorff and étale. 
    By Kumjian--Renault theory (\cite{kumjian1986c,renault2008cartan}), this is equivalent to the existence of a \emph{Cartan subalgebra} $D \subseteq A$, which corresponds to the unit space $\G^{(0)} \subseteq \G$. 

    Two significant questions are whether a groupoid model $\G$ with trivial twist $\Sigma = \G \times \mathbb{T}$ can always be found (\cite[Problem~XLIII]{STW25_99problems}) and whether the groupoid $\G$ can always be taken to be principal, meaning that the only groupoid elements whose range and source coincide are the elements of the unit space $\G^{(0)}$. In this paper, we provide a simultaneous positive answer to both questions in the purely infinite setting. 

    \begin{thm}\label{Thm:Kirchberg}
        Let $A$ be a Kirchberg algebra satisfying the UCT. Then there exists a principal, second-countable, locally compact, Hausdorff, étale groupoid $\G$ with a finite-dimensional unit space such that $A \cong C_r^*(\G)$.
    \end{thm}

    Under Kumjian--Renault theory, principality corresponds to the property that the Cartan subalgebra satisfies the unique extension property for pure states, and so is a \emph{C$^*$-diagonal} in the sense of Kumjian (\cite{kumjian1986c}); see Section~\ref{sec:preliminaries} for definitions. Our motivation for seeking principal groupoid models is that Cartan inclusions without the unique extension property have infinite diagonal dimension (\cite{LLW23}), making them unlikely candidates for classification by $K$-theoretic invariants.
    
    Combining Theorem~\ref{Thm:Kirchberg} with Kumjian--Renault theory and Xin Li's construction in the stably finite setting (\cite{li2020every}), we settle the existence question for C$^*$-diagonals in classifiable C$^*$-algebras. 
    
    \begin{thm}\label{Thm:every-classifiable-diagonal}
        Every simple, separable, nuclear, UCT, Toms--Winter regular\footnote{This means any one of finite nuclear dimension (\cite{WZ10}), Jiang--Su stability (\cite{JS99}), or strict comparison (see for example \cite{GardellaPerera26}). These regularity conditions are now known to be equivalent if $A$ is simple, separable, nuclear and non-elementary (\cite{TomsWinter, Thiel26,OAI:Toms-Winter}). Note: Jiang--Su stability automatically implies $A$ is non-elementary and finite nuclear dimension automatically implies $A$ is nuclear. Since matrix algebras $M_n(\C)$ (resp.\ the compact operators $\K$) have a C$^*$-diagonal consisting of diagonal matrices (resp.\ diagonal compact operators), we have omitted the adjective non-elementary in this theorem.} C$^*$-algebra $A$ has a C$^*$-diagonal $D \subseteq A$. Moreover, $D$ can be chosen to have finite-dimensional spectrum.
    \end{thm}

    Applying Kumjian--Renault theory, every classifiable C$^*$-algebra has a twisted groupoid model that is principal, second-countable, locally compact, Hausdorff and étale.
    In the purely infinite case, the twist can be taken to be trivial by Theorem~\ref{Thm:Kirchberg}.
    Moreover, Xin Li's construction produces a groupoid model with a trivial twist in the stably finite case whenever the $K_0$-group of the C$^*$-algebra is torsion-free \cite[Corollary~1.8]{li2020every}. 

    Theorem~\ref{Thm:every-classifiable-diagonal} also strengthens Xin Li's characterisation of the UCT by the existence of Cartan subalgebras (\cite[Corollary~1.5]{li2020every}) as we can now show a Cartan subalgebra can always be chosen to have the unique extension property for pure states.
    
    \begin{cor}\label{Cor:UCT-diagonal}
        Let $A$ be a simple, separable, nuclear, Toms--Winter regular C$^*$-algebra. Then the following are equivalent:
        \begin{enumerate}[(i)]
            \item\label{it:UCT} $A$ satisfies the UCT;
            \item\label{it:Diagonal} $A$ has a C$^*$-diagonal;
            \item\label{it:Cartan} $A$ has a Cartan subalgebra.
        \end{enumerate}
    \end{cor}

  The implication \eqref{it:UCT}$\Rightarrow$\eqref{it:Diagonal} follows from Theorem~\ref{Thm:every-classifiable-diagonal}, noting that elementary algebras have standard C$^*$-diagonals. The implication \eqref{it:Diagonal}$\Rightarrow$\eqref{it:Cartan} is true by definition, and the implication \eqref{it:Cartan}$\Rightarrow$\eqref{it:UCT} follows from the extension of Tu's theorem (\cite{Tu99}) to twisted groupoids by Barlak--Li (\cite[Theorem~1.1]{barlak2017cartan}) because the existence of a Cartan subalgebra yields an amenable twisted groupoid model for $A$ by Kumjian--Renault theory; see \cite[Corollary~1.2]{barlak2017cartan}.

    The results of this paper represent the culmination of a large body of work constructing principal groupoid models of UCT Kirchberg algebras in various special cases (\cite{hjelmborgstability,rordam2012purely,elliott2016k, suzuki2017amenable,brown2016purely, sibbel2024cantor, ES26, ES26preprint, ES26Stable, LLW23}). The salient point is that the traditional approaches to constructing groupoid models for Kirchberg algebras using graphs and more general graph-like objects, such as higher-rank graphs (\cite{spielberg2007graph,li2019cartan, CFH20}) and topological graphs (\cite{katsura2008class,yeend2006topological}), typically lead to groupoid models that are not principal owing to the presence of cycles in the graph-like object.

    The solution to the problem of cycles is to introduce a minimal dynamical system into the construction --- a technique we shall call \emph{liberation} as it makes the natural dynamics on paths in a topological graph \emph{free}, in the sense of having no non-trivial fixed points. This results in a boundary path groupoid that is \emph{principal}. 
    
    Liberation has its roots in the work of Winter and the second-named author (\cite{sibbel2024cantor}) with the groupoid-theoretic viewpoint brought to the fore in the authors' paper on principal groupoid models of Cuntz algebras (\cite{ES26}). 
    The stable UCT Kirchberg case of Theorem~\ref{Thm:Kirchberg} was proved in \cite{ES26Stable} by the authors with the key insight being to replace the Cantor minimal system used in \cite{sibbel2024cantor,ES26} with a minimal dynamical system, constructed by Deeley--Putnam--Strung (\cite{deeley2018constructing}), which takes place on an infinite compact metric space that is non-contractible but nonetheless has the same $K$-theory as the one-point space.

    In this paper, we prove the unital case of Theorem~\ref{Thm:Kirchberg}, which also gives a second proof of the stable case via taking the product of the groupoid and the full equivalence relation on $\N$. Since every Kirchberg algebra is either unital or stable (\cite{Zhang92}), these two results together prove Theorem~\ref{Thm:Kirchberg}.

    In addition to the aforementioned dynamical system of Deeley--Putnam--Strung, the main ingredients for the proof of Theorem~\ref{Thm:Kirchberg} are Katsura's sufficient conditions for when a topological graph C$^*$-algebra is a UCT Kirchberg algebra and his construction of unital topological graph C$^*$-algebras realising all pairs of countable abelian groups as $K$-groups, including full control over the $K_0$-class of the unit (\cite{katsura2008class}). We prove that liberation using the Deeley--Putnam--Strung dynamical system has no effect on the $K$-theory in general (see Proposition~\ref{prop:product-k-theory}) and that, for the specific topological graph C$^*$-algebras constructed by Katsura, the properties of minimality and contraction are also preserved, which ensures that the corresponding C$^*$-algebra is UCT Kirchberg (see Propositions~\ref{prop:minimality} and \ref{prop:contracting}). 
    The proof of Theorem~\ref{Thm:Kirchberg} then follows by the Kirchberg--Phillips theorem (\cite{Ki95,Ph00}).

    A simpler version of this construction was possible in the stable case (\cite{ES26Stable}), 
    as we could work with a much simpler class of topological graph C$^*$-algebras, where the $K$-theory of the topological graph C$^*$-algebra coincides with that of a locally compact space used in the construction. However, this approach is only valid under $K$-theoretic restrictions in the unital case. The main technical achievement of this paper is the in-depth understanding of the interaction between liberation and Katsura's computation of the $K$-theory of a topological graph C$^*$-algebra using a Toeplitz extension to derive a 6-term exact sequence in $K$-theory.
    
    Interestingly, it is precisely the presence of loops at each vertex in Katsura's construction of topological graphs with prescribed $K$-theory that allows us to prove that the liberated topological graph remains minimal and contracting, which implies that its C$^*$-algebra is a UCT Kirchberg algebra. At the same time, it is precisely this feature that ensures that the un-liberated topological graph has a boundary path groupoid that is not principal!  

    \subsection*{Structure of paper}
    In Section~\ref{sec:preliminaries}, we fix notation and recall some basic results on Kirchberg algebras, groupoids and Kumjian--Renault theory. Section~\ref{sec:topological-graphs} is a review of the theory of topological graphs and their C$^*$-algebras. In Section~\ref{sec:dynamical-modification}, we introduce our liberation technique for topological graphs, and we establish principality of the associated boundary path groupoid. Section~\ref{sec:k-theory} contains the proof that liberation with respect to the Deeley--Putnam--Strung dynamical system preserves $K$-theory. In Section~\ref{sec:principal-models}, we apply liberation to Katsura's topological graph models for Kirchberg algebras. Finally, Section~\ref{sec:main-proofs} contains proofs of our main results.
    
    \subsection*{Acknowledgments}
    The authors would like to thank Wilhelm Winter for his helpful comments after our seminar talks on this project at the University of Münster. 

    \subsection*{AI Statement}
    No AI was used in the development of the mathematics in this paper. The authors have made use of standard AI tools (e.g.\ ChatGPT, Google Gemini) during the write-up phase to accelerate the proofreading and polishing of the manuscript.

	\section{Preliminaries} \label{sec:preliminaries}
	\renewcommand{\thethm}{\arabic{thm}}
	\numberwithin{thm}{section}
    In this section, we fix the notation and recall the unital classification theorem for UCT Kirchberg algebras and the Kumjian--Renault correspondence.

    A \emph{Kirchberg algebra} is a simple, separable, nuclear, purely infinite C$^*$-algebra. The unital Kirchberg--Phillips theorem (\cite{Ki95,Ph00}) states that two unital UCT Kirchberg algebras $A$ and $B$ are isomorphic if and only if $(K_0(A),[1_A]_0,K_1(A)) \cong (K_0(B),[1_B]_0,K_1(B))$. Here, the isomorphism of $K_0$-groups is required to send $[1_A]_0$ to $[1_B]_0$ (see \cite[Theorem~8.4.1]{rordam2002classification}).

    We assume familiarity with the basic theory of topological groupoids and their C$^*$-algebras (\cite{renault1980groupoid,renault2008cartan}). Throughout the paper, groupoids are locally compact and Hausdorff. For a groupoid $\G$, we write $\G^{(0)}$ for its unit space and $r,s:\G\to\G^{(0)}$ for the range and source maps. The groupoid is said to be \emph{étale} if these maps are local homeomorphisms. For $x\in\G^{(0)}$, the \emph{isotropy group} is $\G_x^x=\{g\in\G:r(g)=s(g)=x\}$.
    We call $\G$ \emph{principal} if $\G_x^x=\{x\}$ for every unit $x$, and \emph{topologically principal} if the units with trivial isotropy are dense in $\G^{(0)}$. For an étale groupoid $\G$, the reduced C$^*$-algebra is denoted by $C_r^*(\G)$. More generally, if $\Sigma$ is a twist over $\G$, we write $C_r^*(\G,\Sigma)$ for the associated reduced twisted C$^*$-algebra (see \cite[Section~4]{renault2008cartan}). In either case, we identify $C_0(\G^{(0)})$ with its canonical abelian subalgebra of $C_r^*(\G,\Sigma)$.

    Let $A$ be a C$^*$-algebra. A \emph{Cartan subalgebra} $D\subseteq A$ is a non-degenerate maximal abelian sub-C$^*$-algebra that admits a faithful conditional expectation $A\to D$ and whose normalisers $\mathcal N_A(D) = \{n\in A:n^*Dn\subseteq D \text{ and } nDn^*\subseteq D\}$ generate $A$ as a C$^*$-algebra (\cite[Definition~5.1]{renault2008cartan}).
    A Cartan subalgebra is a \emph{C$^*$-diagonal} if, in addition, every pure state on $D$ has a unique extension to a state on $A$ (see \cite{kumjian1986c}).

    Kumjian--Renault theory identifies Cartan inclusions in separable C$^*$-algebras with twisted, topologically principal, second-countable, locally compact, Hausdorff, étale groupoids (\cite{kumjian1986c,renault2008cartan}).
    More precisely, every such Cartan inclusion is isomorphic to $C_0(\G^{(0)})\subseteq C_r^*(\G,\Sigma)$, where the twisted groupoid $(\G,\Sigma)$ is uniquely determined up to isomorphism by the inclusion.
    Conversely, every twisted groupoid with these properties gives a Cartan inclusion. The inclusion is a C$^*$-diagonal precisely when $\G$ is principal.

	\section{Topological graphs and their \cstar-algebras} \label{sec:topological-graphs}

    In this section, we recall the definition of topological graph C$^*$-algebras via correspondences and introduce their Toeplitz
    algebras. We also describe their boundary path groupoid models, and recall Katsura's sufficient conditions for these algebras to be UCT Kirchberg algebras. Our main references are \cite{katsura2004class,katsura2006class2,katsura2008class,kumjian2017twisted}.

    \subsection{Graphs, paths and boundary paths} \label{subsec:graphs-and-paths}

    A \emph{topological graph} is a 4-tuple  $E=(E^0,E^1,d_E,r_E)$ consisting of two locally compact Hausdorff spaces $E^0$ and $E^1$ together with maps $d_E,r_E: E^1\to E^0$, where $d_E$ is a local homeomorphism, and $r_E$ is continuous (\cite[Definition~2.1]{katsura2004class}). We call $E$ \emph{second-countable} if both $E^0$ and $E^1$ are second-countable.

    A vertex $v\in E^0$ is said to be \emph{regular} if it has a neighbourhood $V$ such that $r_E^{-1}(V)$ is compact and $r_E(r_E^{-1}(V))=V$. The regular vertices form an open subset $E^0_{rg}\subseteq E^0$. Its complement $E^0_{sg}=E^0\setminus E^0_{rg}$ is the set of \emph{singular} vertices.

    For $n\geq 2$, the space of paths of length $n$ is
    \begin{equation}
        E^n=\{(\mu_1,\ldots,\mu_n)\in(E^1)^n: d_E(\mu_i)=r_E(\mu_{i+1}) \text{ for } 1\leq i<n\},
    \end{equation}
    with the subspace topology.
    We regard vertices as paths of length zero and write $E^*=\bigsqcup_{n=0}^{\infty}E^n$ for the space of finite paths, equipped with the disjoint union topology. The infinite path space is
    \begin{equation}\label{eqn:infinite-path-space}
        E^\infty=\{(\mu_i)_{i=1}^{\infty}\in(E^1)^{\mathbb N}: d_E(\mu_i)=r_E(\mu_{i+1}) \text{ for every } i\in\mathbb N\}.
    \end{equation}

    For a path $\mu=(\mu_1,\ldots,\mu_n)$ of finite length, set $r_E(\mu)=r_E(\mu_1)$ and $d_E(\mu)=d_E(\mu_n)$. On $E^0$, both maps are the identity. This defines continuous maps $d_E,r_E: E^*\to E^0$. 
    We also define a continuous map $r_E: E^\infty \to E^0$ by setting $r_E(\mu)=r_E(\mu_1)$. The length of a finite or infinite path is denoted by $|\mu|\in\mathbb N_0\cup\{\infty\}$. For each finite $n$, the map $d_E: E^n\to E^0$ is a local homeomorphism.

    Following \cite[Definition~3.7]{kumjian2017twisted}, the \emph{boundary path space} of $E$ is the set $\partial E= E^\infty \sqcup \{\mu\in E^*:d_E(\mu)\in E^0_{sg}\}$ equipped with the topology defined as follows: For $S\subseteq E^*$, define
    \begin{equation}
        Z(S) = \{\mu\in\partial E: r_E(\mu)\in S \text{ or } (\mu_1,\ldots,\mu_i)\in S \text{ for some } 1\leq i\leq|\mu|\}.
    \end{equation}
    Then the sets $Z(U)\setminus Z(K)$, where $U\subseteq E^*$ is open and $K\subseteq E^* $ is compact, form a basis for a locally compact, Hausdorff topology on $\partial E$.
    If $E$ is second-countable, then so is $\partial E$. In that case its topology is determined by the following description of convergent sequences.

    \begin{lem}[{\cite[Lemma~3.8]{kumjian2017twisted}}] \label{lem:Kumjian-Li}
        Let $E$ be a topological graph, let $(\mu^{(n)})_{n=1}^{\infty}$ be a sequence in $\partial E$, and let $\mu\in\partial E$. Then $\mu^{(n)}\to\mu$ if and only if the following hold:
        \begin{enumerate}[(i)]
            \item $r_E(\mu^{(n)})\to r_E(\mu)$;
            \item for every integer $i$ with $1\leq i\leq|\mu|$, there exists $N \ge 1$ such that $|\mu^{(n)}| \ge i$ for all $n \ge N$ and
            $(\mu^{(n)}_1,\ldots,\mu^{(n)}_i) \longrightarrow (\mu_1,\ldots,\mu_i)$;
         
            \item if $|\mu|<\infty$, then the set
            $\{n:|\mu^{(n)}|>|\mu| \text{ and } \mu^{(n)}_{|\mu|+1}\in K\}$ is finite for every compact subset $K\subseteq E^1$.
        \end{enumerate}
    \end{lem}

    \subsection{Topological graph \cstar-algebras} \label{subsec:graph-correspondences}

    Let $E=(E^0,E^1,d_E,r_E)$ be a topological graph. We recall the correspondence associated to $E$ and the definitions of its Toeplitz and graph C$^*$-algebras, following \cite{katsura2004class} and \cite{katsura2008class}.

    Define $C_d(E^1) = \{\xi\in C(E^1): \langle\xi,\xi\rangle\in C_0(E^0)\}$,
    where
    \begin{equation}
        \langle\xi,\eta\rangle(v) = \sum_{e\in d_E^{-1}(v)} \overline{\xi(e)}\,\eta(e)
    \end{equation}
    for $\xi, \eta \in C_d(E^1) $ and $v \in E^0.$ 
    The left and right actions of $C_0(E^0)$ on $C_d(E^1)$ are given by
    \begin{equation}
    (f\xi g)(e)= f(r_E(e))\xi(e)g(d_E(e)),
    \end{equation}
    for $f,g\in C_0(E^0)$, $\xi\in C_d(E^1)$ and $e\in E^1$.
    The right multiplication and inner product make $C_d(E^1)$ a Hilbert $C_0(E^0)$-module; see \cite[Proposition~1.10]{katsura2004class}. Moreover, $C_c(E^1)$ is contained in and dense in $C_d(E^1)$; see \cite[Lemmas~1.5 and~1.6]{katsura2004class}. We note that $\|\xi\|=\|\langle\xi,\xi\rangle\|^{1/2}_{C_0(E^0)}$ is the norm on $C_d(E^1)$.

    For a Hilbert module $Y$, we write $\mathcal L(Y)$ for the C$^*$-algebra of adjointable operators on $Y$ and $\K(Y)$ for the closed linear span of the rank-one operators given by $\theta_{\xi,\eta}(\zeta)=\xi\langle\eta,\zeta\rangle$ for $\xi,\eta,\zeta \in Y$.
    The left action of $C_0(E^0)$ on $C_d(E^1)$ defines a $^*$-homomorphism $\pi_{r_E}: C_0(E^0)\longrightarrow\mathcal L(C_d(E^1))$,
    whose restriction to $C_0(E^0_{rg})$ is an injective $^*$-homomorphism into $\K(C_d(E^1))$.

    A \emph{Toeplitz $E$-pair} $(T_E^0,T_E^1)$ in a C$^*$-algebra $A$ consists of a $^*$-homomorphism $T_E^0: C_0(E^0)\to A$ and a linear map $T_E^1: C_d(E^1)\to A$ satisfying
    \begin{align}
        T_E^1(\xi)^*T_E^1(\eta) &= T_E^0(\langle\xi,\eta\rangle),\\
        T_E^0(f)T_E^1(\xi) &= T_E^1(\pi_{r_E}(f)\,\xi),
    \end{align}
    for $f\in C_0(E^0)$ and $\xi,\eta\in C_d(E^1)$; see \cite[Definition~2.2]{katsura2004class}. 
    Such a pair induces a $^*$-homomorphism $\Phi_E:\K(C_d(E^1))\to A$ via
    \begin{equation}
         \Phi_E(\theta_{\xi,\eta}) = T_E^1(\xi)\,T_E^1(\eta)^*.
    \end{equation}

    The \emph{Toeplitz algebra} $\mathcal T(E)$ is the universal C$^*$-algebra generated by a Toeplitz $E$-pair. A Toeplitz $E$-pair $T_E=(T_E^0,T_E^1)$ is called a \emph{Cuntz--Krieger $E$-pair} if $T_E^0(f)=\Phi_E(\pi_{r_E}(f))$ for all $f\in C_0(E^0_{rg})$. 

    \begin{dfn}[{\cite[Definition~2.10]{katsura2004class}}] \label{dfn:TopologicalGraphAlgebra}
        The \emph{topological graph C$^*$-algebra} $\mathcal O(E)$ is the universal C$^*$-algebra generated by a Cuntz--Krieger $E$-pair $t_E=(t_E^0,t_E^1)$.
    \end{dfn}

  The coefficient map $t_E^0: C_0(E^0)\to\mathcal O(E)$ is injective by \cite[Proposition~3.7]{katsura2004class}. Moreover, the topological graph C$^*$-algebra $\mathcal O(E)$ is unital whenever $E^0$ is compact and in this case $1_{\mathcal O(E)}=t_E^0(1_{C(E^0)})$;
    see \cite[Proposition~7.1]{katsura2006class2}.
    Finally, since every Cuntz--Krieger $E$-pair is a Toeplitz $E$-pair, there is a canonical surjective $^*$-homomorphism $q_E:\mathcal T(E)\rightarrow\mathcal O(E)$ with $q_E \circ T_E^i=t_E^i$ for $i\in\{0,1\}.$
    We will discuss the kernel of this map and the resulting Toeplitz extension further in Subsection~\ref{subsec:fock}.

    \subsection{Boundary path groupoids} \label{subsec:top-graphs-groupoids}
    Let $E=(E^0,E^1,d_E,r_E)$ be a second-countable topological graph. The boundary path space carries a natural shift map $\sigma_E$. Its domain is $\dom(\sigma_E)=\partial E\setminus E^0_{sg},$
    and it is defined by
    \begin{equation}
    \sigma_E(\mu) = 
        \begin{cases}
            d_E(\mu_1), & |\mu|=1,\\
            (\mu_2,\mu_3,\ldots), & |\mu|\geq 2.
        \end{cases}
    \end{equation}
    By \cite[Lemma~6.1]{kumjian2017twisted}, $\sigma_E$ is a local homeomorphism from its open domain onto an open subset of $\partial E$.

    We recall the groupoid construction associated to such a map introduced by Deaconu (\cite{deaconu1995groupoids}) and Renault (\cite{renault2000cuntz}).
    Let $X$ be a second-countable, locally compact, Hausdorff space and let $\sigma:\dom(\sigma) \to\ran(\sigma)$ be a local homeomorphism between open subsets of $X$. 
    For $n\geq1$, the map $\sigma^n$ is the $n$-fold composition of $\sigma$, defined on the open subset of $X$ where this makes sense, and $\sigma^0=\id_X$. 
    The \emph{Deaconu--Renault groupoid} $\G(X,\sigma)$ then is
    \begin{equation}
         \{(x,n-m,y):
        n,m\in\mathbb N_0,\,
        x\in\dom(\sigma^n),\,
        y\in\dom(\sigma^m),\,
        \sigma^n(x)=\sigma^m(y)\}.
    \end{equation}
    We identify its unit space with $X$ via $x\mapsto(x,0,x)$. The range and source maps are $r(x,k,y)=x$ and $s(x,k,y)=y$,
    and the groupoid operations are given by
    \begin{equation}
        (x,k,y)\,(y,\ell,z)=(x,k+\ell,z), \qquad (x,k,y)^{-1}=(y,-k,x).
    \end{equation}
    A basis for the topology consists of the sets
    \begin{equation}
        B(U,n,m,V) = \{(x,n-m,y):
        x\in U,\,
        y\in V,\,
        \sigma^n(x)=\sigma^m(y)\},
    \end{equation}
    where $U\subseteq\dom(\sigma^n)$ and
    $V\subseteq\dom(\sigma^m)$ are open and the restrictions $\sigma^n|_U$ and $\sigma^m|_V$ are injective. This makes $\G(X,\sigma)$ a second-countable, locally compact, Hausdorff, étale groupoid for which the full and reduced C$^*$-algebras coincide.

    Thus, for a second-countable topological graph $E$, we can consider the associated Deaconu--Renault groupoid $\G(\partial E,\sigma_E)$ which is called the \emph{boundary path groupoid} and gives a model for the corresponding topological graph algebra.

    \begin{prop}[{\cite[Theorem~5.2]{yeend2006topological}; \cite[Theorem~6.7]{kumjian2017twisted}}] \label{prop:groupoid-model-for-graphs}
        Let $E$ be a second-countable topological graph. Then
        \begin{equation}
            \mathcal O(E)\cong C_r^* \bigl( \G(\partial E,\sigma_E) \bigr).
        \end{equation}
    \end{prop}

    We will be particularly interested in the principality of the groupoid $\G(\partial E,\sigma_E)$. We record the following simple observation here as a lemma. 

    \begin{lem}\label{lem:DR-principality}
        Let $E$ be a second-countable topological graph. Then $\G(\partial E,\sigma_E)$ is principal if and only if, for all $n,m \in \N_0$ and all infinite paths $\mu \in  E^\infty$, we have that $\sigma_E^n(\mu)=\sigma_E^m(\mu)$ already implies $m=n$.
    \end{lem}
    \begin{proof}
        The groupoid $\G(\partial E,\sigma_E)$ is principal if and only if, for every $\mu\in \partial E$ and $n,m\in\N_0$ for which both $\sigma_E^n(\mu)$ and $\sigma_E^m(\mu)$ are defined, we have that $\sigma_E^n(\mu)=\sigma_E^m(\mu)$ implies $m=n$.
        However, this implication holds automatically on finite paths as each application of the shift decreases the path length by one. Therefore, we only need to consider infinite paths, as claimed.
    \end{proof}

    \subsection{Minimality and contraction} \label{subsec:minimality-contraction}
    Let $E=(E^0,E^1,d_E,r_E)$ be a topological graph.
    We recall Katsura's sufficient conditions for a topological graph C$^*$-algebra to be a UCT Kirchberg algebra (\cite{katsura2008class}).\footnote{Note that Katsura's definition of a Kirchberg algebra already includes the UCT.}

    For $v\in E^0$, its \emph{positive orbit} is
    \begin{equation}
        \operatorname{Orb}^+(v) = \{r_E(\mu):\mu\in E^*,\, d_E(\mu)=v\}.
    \end{equation}
    A \emph{negative orbit} of $v$ is a boundary path $\mu\in\partial E$ with $r_E(\mu)=v$. For such a path, set
    \begin{equation}
        \operatorname{Orb}^-(v,\mu) = \{v\} \cup \{d_E(\mu_i):i \in \N,\, 1\leq i\leq|\mu|\},
    \end{equation}
    and define the orbit space of $v$ with respect to a negative orbit $\mu$ of $v$ by 
    \begin{equation}
        \operatorname{Orb}(v,\mu) = \bigcup_{w\in\operatorname{Orb}^-(v,\mu)} \operatorname{Orb}^+(w).
    \end{equation}

    \begin{dfn}[{\cite[Definition~1.8]{katsura2008class}}] \label{dfn:graph-minimality}
        A topological graph $E$ is \emph{minimal} if $\operatorname{Orb}(v,\mu)$ is dense in $E^0$ for every $v\in E^0$ and every negative orbit $\mu$ of $v$.
    \end{dfn}
    
    To formulate the contraction condition, let $U\subseteq E^n$ and $U'\subseteq E^m$ for $n,m\geq 1$. For $k=\min\{n,m\}$, we then write $U|_k = \{(\mu_1,\ldots,\mu_k):\mu\in U\}$ and 
    \begin{equation}
        U\pitchfork U'=(U|_k)\cap(U'|_k).
    \end{equation}

    \begin{dfn}[{\cite[Definitions~2.3 and~2.7]{katsura2008class}}] \label{dfn:graph-contraction}
        Let $E = (E^0,E^1,d_E,r_E)$ be a topological graph. A non-empty open subset $V\subseteq E^0$ is called \emph{contracting} if $\overline V$ is compact and there exist $N \in \N$ and non-empty open sets $U_k\subseteq E^{n_k}$ for $ 1\leq k\leq N$, where $n_1,\ldots,n_N$ are positive integers, such that
        \begin{enumerate}[(i)]
            \item $r_E(U_k)\subseteq V$ for $1\leq k\leq N$;
            \item   $U_k\pitchfork U_l=\emptyset$ whenever $k\neq l$;
            \item $\overline V\subsetneq\bigcup_{k=1}^N d_E(U_k)$.
        \end{enumerate}
        The graph $E$ is said to be \emph{contracting at $v_0\in E^0$} if $\overline{\operatorname{Orb}^+(v_0)}=E^0$ and every neighbourhood of $v_0$ contains a contracting open set.
        We call $E$ \emph{contracting} if it is contracting at some vertex.
    \end{dfn}

    The following result combines these conditions.

    \begin{thm}[{\cite[Corollary~B]{katsura2008class}}]
    \label{thm:graph-contracting-Kirchberg-algebra}
        Let $E$ be a second-countable, minimal, contracting topological graph.
        Then $\mathcal O(E)$ is a UCT Kirchberg algebra.
    \end{thm}

    \section{Liberating topological graphs} \label{sec:dynamical-modification}

    In this section, we develop the general theory of liberation for topological graphs, which we will use to prove the main results of this paper. Special cases of this general construction were already used by the authors in \cite{ES26, ES26preprint, ES26Stable} and the key idea can be traced back further to work of Winter and the second-named author in \cite{sibbel2024cantor}.

    We begin with the definition of the modified topological graph with respect to a dynamical system, which we state for a general homeomorphism of a locally compact Hausdorff space. 
    \begin{dfn}\label{def:modified-topological-graph} 
    Let $E=(E^0,E^1,d_E,r_E)$ be a topological graph, let $X$ be a locally compact Hausdorff space, and let $\rho: X\to X$ be a homeomorphism.
    We define the topological graph
    \begin{equation}
        E\times_\rho X = ((E\times_\rho X)^0,(E\times_\rho X)^1,d,r)
    \end{equation}
    by setting $(E\times_\rho X)^0 = E^0\times X$, $(E\times_\rho X)^1 = E^1\times X$, $d(e,x) = (d_E(e),x)$ and $r(e,x) = (r_E(e),\rho(x))$.
    \end{dfn}
    We note straight away that $E\times_\rho X$ is a valid topological graph: since $\rho$ is a homeomorphism, $d_E$ is a local homeomorphism, and $r_E$ is continuous, it follows that $d$ is a local homeomorphism and $r$ is continuous.
    Moreover, if $E$ and $X$ are second-countable, then so is $E\times_\rho X$. 

    In the following lemma, we describe the regular vertices and singular vertices of the modified graph when the dynamics takes place on a compact space.
    \begin{lem}\label{lem:product-regular-vertices}
        Let $E$ be a topological graph and let $\rho: X\to X$ be a homeomorphism of a compact Hausdorff space $X$. Then $(E\times_\rho X)^0_{rg}=E^0_{rg}\times X$ and $(E\times_\rho X)^0_{sg}=E^0_{sg}\times X$.
    \end{lem}

    \begin{proof}
        Let $F=E\times_\rho X$. First, let $(v,x)\in E^0_{rg}\times X$. Then, there exists a neighbourhood $V\subseteq E^0$ of $v$ such that $r_E^{-1}(V)$ is compact and $r_E(r_E^{-1}(V))=V$. Now $V\times X$ is a neighbourhood of $(v,x)$, and $r_F^{-1}(V\times X)=r_E^{-1}(V)\times X$ is compact. Moreover, $r_F(r_F^{-1}(V\times X)) = r_E(r_E^{-1}(V))\times\rho(X) = V\times X.$
        Thus $(v,x)\in F^0_{rg}$.

        Conversely, suppose that $(v,x)\in F^0_{rg}$. Then, there exists a neighbourhood $W\subseteq E^0\times X$ of $(v,x)$ such that $r_F^{-1}(W)$ is compact and $r_F(r_F^{-1}(W))=W$. Choose closed neighbourhoods $V\subseteq E^0$ of $v$ and $U\subseteq X$ of $x$ such that $(v,x)\in V\times U\subseteq W$.
        Since $V\times U$ is closed, the set $r_F^{-1}(V\times U) = r_E^{-1}(V)\times\rho^{-1}(U)$ is a closed subset of $r_F^{-1}(W)$ and hence compact. As $\rho^{-1}(U)$ is non-empty, projection onto the first coordinate shows that $r_E^{-1}(V)$ is compact.
        Finally, let $v'\in V$.
        Since $(v',x)\in W=r_F(r_F^{-1}(W))$, there exists $(e,x')\in F^1$ such that $r_F(e,x')=(v',x).$
        In particular, $r_E(e)=v'$. Thus $r_E(r_E^{-1}(V))=V$, and $v\in E^0_{rg}$. 
        
        Therefore, we have shown that $F^0_{rg}=E^0_{rg}\times X$.
        Taking complements now gives that $F^0_{sg}=E^0_{sg}\times X$.
    \end{proof}

    Next, we identify the paths in the modified topological graph. We shall call the homeomorphisms $\Phi_n:E^n\times X\longrightarrow (E\times_\rho X)^n$ defined in the following lemma the \emph{path space homeomorphisms}. 
    For completeness, we define $\Phi_0$ to be the identity map $E^0\times X\to (E \times_\rho X)^0$.
    
    \begin{lem}\label{lem:path-homeo}
    Let $E$ be a topological graph, let $X$ be a locally compact Hausdorff space, and let $\rho: X\to X$ be a homeomorphism.
    Let $n\geq 1$. Define $\Phi_n: E^n\times X\longrightarrow (E\times_\rho X)^n$ by
    \begin{equation}
        \Phi_n((e_1,\ldots,e_n),x) = \bigl(
        (e_1,\rho^{-1}(x)),
        (e_2,\rho^{-2}(x)), \ldots,
        (e_n,\rho^{-n}(x))
        \bigr),
    \end{equation}
    where $E^n$ is the space of paths of length $n$ in $E$ and $(E\times_\rho X)^n$ is the space of paths of length $n$ in $(E\times_\rho X)$.
    Then $\Phi_n$ is a homeomorphism.
    Moreover, we have
    \begin{equation}\label{eqn:modified-path-domain-range}
        r(\Phi_n(\mu,x))=(r_E(\mu),x)
        \quad \text{and} \quad
        d(\Phi_n(\mu,x))=(d_E(\mu),\rho^{-n}(x)).
    \end{equation}
    \end{lem}
    \begin{proof}
    Since $\rho$ is a homeomorphism, the map $\Phi_n$ is clearly injective and continuous. The path relation in $E\times_\rho X$ implies that $\Phi_n$ is surjective. Indeed, the inverse $\Phi_n^{-1}$ sends $((e_1,x_1),\ldots,(e_n,x_n))$ to $((e_1,\ldots,e_n),\rho(x_1))$. Since $\Phi_n^{-1}$ is continuous, it follows that $\Phi_n$ is a homeomorphism. It follows from the definition of $\Phi_n$ that \eqref{eqn:modified-path-domain-range} holds.
    \end{proof}

    Now we get to the key point of the liberation technique. We show that the absence of periodic points for $\rho$ ensures principality of the boundary path groupoid of $E\times_\rho X$.
    \begin{prop}\label{prop:principal-groupoid}
        Let $E$ be a second-countable topological graph, and let $\rho: X\to X$ be a homeomorphism of a second-countable, locally compact, Hausdorff space $X$.
        Suppose that $\rho$ has no periodic points. Set $F=E\times_\rho X$. Then $\G=\G(\partial F,\sigma_F)$ is a principal, second-countable, locally compact, Hausdorff, étale groupoid, and $\mathcal O(F)\cong C_r^*(\G).$
    \end{prop}

    \begin{proof}
        Since $F$ is second-countable, $\G$ is a second-countable, locally compact, Hausdorff, étale groupoid, and Proposition~\ref{prop:groupoid-model-for-graphs} gives $\mathcal O(F)\cong C_r^*(\G)$.

        By Lemma~\ref{lem:DR-principality}, it remains to show that no infinite path in $F$ has two distinct equal shifts.
        Let $\mu\in F^\infty$. The path compatibility condition \eqref{eqn:infinite-path-space} implies that
        \begin{equation}
            \mu=\bigl(e_i,\rho^{-i}(x)\bigr)_{i=1}^\infty
        \end{equation}
        for some $(e_i)_{i=1}^\infty\in E^\infty$ and $x\in X$. Note that here, $x$ is the $X$-coordinate of $r_F(\mu)$. If $\sigma_F^n(\mu)= \sigma_F^m(\mu)$ for $n,m\in\N_0$, comparison of the $X$-coordinates of their ranges gives $\rho^{-n}(x)=\rho^{-m}(x).$
        Since $\rho$ has no periodic points, we obtain $n=m$. Thus, $\G$ is principal.
    \end{proof}

    Recall that a homeomorphism $\rho: X\to X$ is said to be \emph{minimal} if the only closed $\rho$-invariant subsets of $X$ are $\emptyset$ and $X$.
    If $X$ is infinite and Hausdorff, every minimal homeomorphism of $X$ has no periodic points. Indeed, a periodic point would have a finite orbit, which would be a proper non-empty closed $\rho$-invariant set of $X$.
    Thus, Proposition~\ref{prop:principal-groupoid} applies whenever $\rho$ is minimal and $X$ is infinite.

    In the proof of our main results,
    we shall use the dynamical systems constructed by Deeley, Putnam and Strung in their work on a dynamical presentation of the Jiang--Su algebra (\cite{deeley2018constructing}); see also their follow-up paper on minimal dynamics (\cite{deeley2023min-homeo}).
    We record the properties of these dynamical systems needed here.

    \begin{thm}[{\cite[Corollary~1.12]{deeley2018constructing}}]
    \label{thm:DPS-Theorem} 
        There exist an infinite compact metric space $Z$ with covering dimension at most three and a minimal homeomorphism $\rho: Z\to Z$ such that the canonical inclusion $ \iota:\C\longrightarrow C(Z)$ given by $\lambda\longmapsto\lambda 1_{C(Z)}$ is a $KK$-equivalence.
        In particular, $K_0(C(Z))\cong\mathbb Z$ and $ K_1(C(Z))=0$.
    \end{thm}

    We summarise the main results of this section as a proposition.
    \begin{prop}\label{prop:liberation-DPS}
        Let $(Z,\rho)$ be a Deeley--Putnam--Strung minimal dynamical system as in Theorem~\ref{thm:DPS-Theorem}.
        Let $E$ be a second-countable topological graph.
        Let $\mathcal{G}$ denote the boundary path groupoid of $E\times_\rho Z$. Then $\mathcal{G}$ is a principal, second-countable, locally compact, Hausdorff, étale groupoid with $\mathcal{O}(E\times_\rho Z) \cong C_r^*(\mathcal{G})$.
    \end{prop}
    \begin{proof}
        Since $Z$ is infinite and $\rho$ is minimal, $\rho$ has no periodic points. Hence, Proposition~\ref{prop:principal-groupoid} yields the claim.
    \end{proof}
 
    \section{Preservation of \safeK-theory} \label{sec:k-theory}

    The main goal of this section is to show that the liberation technique introduced in Section~\ref{sec:dynamical-modification} preserves the $K$-theory of the topological graph C$^*$-algebra, including the position of the unit, when the minimal dynamical system constructed by Deeley--Putnam--Strung is used.
    
    \subsection{The induced homomorphisms}  \label{subsec:induced-homomorphisms}

    Given a topological graph $E$ and a homeomorphism $\rho:X \rightarrow X$ on a compact Hausdorff space, there is an obvious projection map $E\times_\rho X \rightarrow E$ given by $(v,x) \mapsto v$ on vertices and $(e,x) \mapsto e$ on edges.
    
    We begin by verifying that this projection is a \emph{factor map} in the sense of Katsura (\cite[Definition~2.1]{katsura2006class2}) and so it induces a $^*$-homomorphism at the level of the Toeplitz algebras $\widetilde\alpha:\T(E) \rightarrow \T(E\times_\rho X)$. We also show that the projection map is \emph{regular} (\cite[Definition~2.6]{katsura2006class2}), so the induced Toeplitz-algebra homomorphism descends to a $^*$-homomorphism $\alpha:\O(E) \rightarrow \O(E\times_\rho X)$.

    For ease of notation, we  make the following identifications of C$^*$-algebras and corresponding Hilbert modules: 
    \begin{align}
        C_0((E\times_\rho X)^0) &\cong C_0(E^0) \otimes C(X),\label{eq:level-0-identifictaion}\\
        C_d((E\times_\rho X)^1) &\cong C_d(E^1) \boxtimes C(X),\label{eq:level-1-identifictaion}
    \end{align}
    where $\boxtimes$ denotes the external tensor product and the $C(X)$ in \eqref{eq:level-1-identifictaion} is the trivial Hilbert-$C(X)$-module (for definitions, see \cite{lance1995hilbert}). These isomorphisms follow from Definition~\ref{def:modified-topological-graph} along with the definition of the C$^*$-correspondences associated to topological graphs (see Subsection~\ref{subsec:graph-correspondences}). 

    \begin{lem}\label{lem:induced-graph-homomorphisms}
        Let $E$ be a topological graph, let $X$ be a compact Hausdorff space, and let $\rho: X\to X$ be a homeomorphism.
        \begin{enumerate}[(i)]
            \item There exists a unique $^*$-homomorphism $\widetilde\alpha:\T(E) \longrightarrow\T(E\times_\rho X)$ such that $\widetilde\alpha(T_E^0(f)) = T_{E\times_\rho X}^0(f \otimes 1_{C(X)})$ and $\widetilde\alpha(T_E^1(\xi)) = T_{E\times_\rho X}^1(\xi \boxtimes  1_{C(X)})$ for all $f \in C_0(E^0)$ and $\xi \in C_d(E^1)$, where $T_E=(T_E^0,T_E^1)$ and $T_{E\times_\rho X}=(T_{E\times_\rho X}^0,T_{E\times_\rho X}^1)$ are the universal Toeplitz pairs. 
            \item There exists a unique $^*$-homomorphism $\alpha:\O(E) \longrightarrow \O(E\times_\rho X)$ such that $\alpha(t_E^0(f)) = t_{E\times_\rho X}^0(f \otimes 1_{C(X)})$ and $\alpha(t_E^1(\xi)) = t_{E\times_\rho X}^1(\xi \boxtimes 1_{C(X)})$ for all $f \in C_0(E^0)$ and $\xi \in C_d(E^1)$, where $t_E=(t_E^0,t_E^1)$ and $t_{E\times_\rho X}=(t_{E\times_\rho X}^0,t_{E\times_\rho X}^1)$ are the universal Cuntz--Krieger pairs.
            \item  If $q_E:\mathcal T(E)\to\mathcal O(E)$ and $q_{E\times_\rho X}:\mathcal T(E\times_\rho X)\to\mathcal O(E\times_\rho X)$ are the canonical quotient maps, then
        $q_{E\times_\rho X}\circ\widetilde\alpha=\alpha\circ q_E$.
            \item If $E^0$ is compact, then $\widetilde\alpha$ and $\alpha$ are unital.
        \end{enumerate}
    \end{lem}
    \begin{proof}
        Set $\Gamma = E\times_\rho X$. Let $m^0:\Gamma^0 \rightarrow E^0$ and $m^1:\Gamma^1 \rightarrow E^1$ be given by $m^0(v,x) = v$ and $m^1(e,x) = e$.
        Note that, under our identifications \eqref{eq:level-0-identifictaion} and \eqref{eq:level-1-identifictaion}, we have $f \circ m^0 = f \otimes 1_{C(X)}$ for $f \in C_0(E^0)$ and $\xi \circ m^1 = \xi \boxtimes 1_{C(X)}$ for all $\xi \in C_d(E^1)$.
        
        Ad (i): We prove that $m=(m^0,m^1)$ is a factor map. It will then follow from \cite[Proposition~2.3]{katsura2006class2} that there exists a unique $^*$-homomorphism  
        $\widetilde\alpha:\T(E) \longrightarrow\T(\Gamma)$ with the required properties.
        Since $X$ is compact, both $m^0$ and $m^1$ are proper. They therefore extend continuously to the respective one-point compactifications by sending $\infty$ to $\infty$. We use the same notation for these extensions.

        We verify the two conditions in \cite[Definition~2.1]{katsura2006class2}. First, for $(e,x)\in\Gamma^1$, we have
        \begin{equation}\begin{split}
            r_E(m^1(e,x)) &=r_E(e) =m^0(r_E(e),\rho(x)) =m^0(r_\Gamma(e,x)), \text{ and}\\
            d_E(m^1(e,x)) &=d_E(e) =m^0(d_E(e),x) =m^0(d_\Gamma(e,x)).
        \end{split}\end{equation}
        Second, let $e\in E^1$ and $(v,x)\in\Gamma^0$ satisfy $d_E(e)=m^0(v,x)=v$. Then $(e,x)$ is the unique edge in $\Gamma^1$ satisfying $m^1(e,x)=e$ and $d_\Gamma(e,x)=(v,x)$.
        Thus $m=(m^0,m^1)$ is a factor map. 

        Ad (ii):  We prove that $m=(m^0,m^1)$ is regular as in \cite[Definition~2.6]{katsura2006class2}. It will then follow from \cite[Proposition~2.9]{katsura2006class2} that 
        there exists a unique $^*$-homomorphism  $\alpha:\O(E) \longrightarrow\O(\Gamma)$ with the required properties.
        Let $(v,x)\in\Gamma^0$ with $m^0(v,x)=v\in E^0_{rg}$.
        Since $v$ is regular, there exists $e\in E^1$ with $r_E(e)=v$. Hence, $ r_\Gamma(e,\rho^{-1}(x))=(v,x),$ so $r_\Gamma^{-1}(v,x)$ is non-empty.
        Moreover, $r_\Gamma^{-1}(v,x)\subseteq\Gamma^1=(m^1)^{-1}(E^1)$ and thus, $m$ is regular.

        Ad (iii): To prove that $q_\Gamma\circ\widetilde\alpha=\alpha\circ q_E$, it suffices to check the identity on the images of the universal Toeplitz pair for $E$. Since $q_E \circ T_E^i = t_E^i$ and $q_\Gamma \circ T_\Gamma^i = t_\Gamma^i$ for $i \in \{0,1\}$, this follows from (i) and (ii).

        Ad (iv): If $E^0$ is compact, then $T_E^0$ and $t_E^0$ are unital. Since $X$ is compact, it follows that $\Gamma^0$ is also compact. Hence, $T_\Gamma^0$ and $t_\Gamma^0$ are unital. 
        Therefore,
        \begin{equation}
            \widetilde\alpha(1_{\T(E)}) =
            \widetilde\alpha(T_E^0(1_{C(E^0)})) =
            T_\Gamma^0(1_{C(E^0)}\otimes 1_{C(X)}) =
            T_\Gamma^0(1_{C(\Gamma^0)}) =
            1_{\T(\Gamma)},
        \end{equation}
        and similarly $\alpha(1_{\O(E)}) = 1_{\O(\Gamma)}$.
    \end{proof}

    Set $I_E = \Ker(q_E)$ and $I_{E \times_\rho X} = \Ker(q_{E \times_\rho X})$.
    By Lemma~\ref{lem:induced-graph-homomorphisms}, we have 
    \begin{equation}
        q_{E \times_\rho X}(\widetilde\alpha(I_E)) = \alpha(q_E(I_E))  = \alpha (0) = 0.
    \end{equation}
    Therefore, $\widetilde\alpha$ restricts to a $^*$-homomorphism $\beta:I_E\longrightarrow I_{E \times_\rho X}$.
    For our $K$-theoretic computation, it will be important to understand the map $\beta$ in more detail. To this end, we will need to consider the Fock space representation of the Toeplitz algebras $\T(E)$ and $\T(E \times_\rho X)$.

\subsection{Fock space representations}\label{subsec:fock}

    Let $E$ be a topological graph.
    For each $n \in \N$, the $n$-fold interior tensor product $C_d(E^1) \otimes_{\pi_{r_E}} \cdots \otimes_{\pi_{r_E}} C_d(E^1)$ of the C$^*$-correspondence $C_d(E^1)$ with respect to the left action $\pi_{r_E}$ will be canonically identified with $C_{d}(E^n)$; see \cite[Section~2]{katsura2004class} and \cite{lance1995hilbert} for definitions. We also write $C_{d}(E^0)$ for the space $C_0(E^0)$ equipped with the trivial Hilbert-$C_0(E^0)$-module structure with left action given by multiplication.
    Taking the direct sum over all $n \in \N_0$, we identify the Fock space $\bigoplus_{n=0}^\infty C_{d}(E^n)$ with $C_d(E^*)$.
    
    Following \cite[Section~3]{katsura2004class}, we can define a Toeplitz pair $(T_E^0,T_E^1)$ in $\L(C_d(E^*))$ as follows: Using left multiplication, we define $T_E^0:C_0(E^0) \rightarrow \L(C_d(E^*))$ for $f \in C_0(E^0)$ by $T_E^0(f)g = fg$ for $g \in C_{d}(E^0)$ and
    \begin{equation}
        T_E^0(f)(\xi \otimes \eta) = (\pi_{r_E}(f)\xi) \otimes \eta,
    \end{equation}
    for all $n \in \N_0$, $\xi \in C_{d}(E^1)$, and $\eta \in C_{d}(E^n)$, and we define $T_E^1:C_d(E^1)\rightarrow \L(C_d(E^*))$ via creation operators, i.e.\
    \begin{equation}
        T_E^1(\xi)(\eta) = \xi \otimes \eta
    \end{equation} 
    for all $\xi \in C_d(E^1)$, $n \in \N_0$, and $\eta \in C_{d}(E^n)$.
    This Toeplitz $E$-pair is universal by \cite[Proposition~3.8]{katsura2004class}, so we may identify $\T(E)$ with the C$^*$-subalgebra of $\L(C_d(E^*))$ generated by the image of $T_E^0$ and $T_E^1$.

    Let $\rho:X \rightarrow X$ be a homeomorphism of a compact Hausdorff space.
    Then the construction above also applies to the topological graph ${E \times_\rho X}$. We write $(T_{E \times_\rho X}^0,T_{E \times_\rho X}^1)$ for the corresponding Toeplitz pair in $\L(C_d({(E \times_\rho X)}^*))$ and identify $\T(E \times_\rho X)$ with the C$^*$-subalgebra of $\L(C_d({(E \times_\rho X)}^*))$ generated by the image of $T_{E \times_\rho X}^0$ and $T_{E \times_\rho X}^1$.

    To keep the notation manageable, we continue to make use of the identifications \eqref{eq:level-0-identifictaion} and \eqref{eq:level-1-identifictaion}. We note that the left action $\pi_{r}$ on the C$^*$-correspondence $C_d((E \times_\rho X)^{1})$ can now be written as $\pi_{r_E} \otimes L_\rho$, where $L_\rho(h):C(X) \rightarrow C(X)$ is left multiplication by $h \circ \rho$ for each $h \in C(X)$; see Definition~\ref{def:modified-topological-graph}.
    
    For $n \in \N$, taking the $n$-fold interior tensor product, we obtain an isomorphism of C$^*$-correspondences $C_{d}((E \times_\rho X)^n) \cong C_{d}(E^n) \boxtimes C(X)$. Under this identification, the left action is given by $\pi_{r_E} \otimes L_{\rho^n}$, where $\pi_{r_E}(f)$ here denotes left multiplication by $f\circ r_E$ on $C_d(E^n)$ for $f\in C_0(E^0)$. 
    Note that these identifications use the $X$-coordinate of the domain of a path.
    
    Hence, we have $C_{d}((E \times_\rho X)^*) \cong C_{d}(E^*) \boxtimes C(X)$. Therefore, we may view 
    $\L(C_{d}(E^*)) \otimes C(X) \subseteq \L(C_{d}((E \times_\rho X)^*))$; see for example \cite{lance1995hilbert}. Since $L_\rho(1_{C(X)}) = \mathrm{id}_{C(X)}$, we have
    \begin{equation}
        T_{E \times_\rho X}^0(f \otimes 1_{C(X)}) =  T_{E}^0(f) \otimes 1_{C(X)}
    \end{equation}
    for all $f \in C_0(E^0)$. Moreover, we have
    \begin{equation}
         T_{E \times_\rho X}^1(\xi \boxtimes 1_{C(X)}) = T_{E}^1(\xi) \otimes 1_{C(X)}
    \end{equation}
    for all $\xi \in C_d(E^1)$ because the creation operator defined on the $n$-fold interior tensor product $C(X) \otimes_{L_\rho} \cdots \otimes_{L_\rho} C(X)$ mapping into the $(n+1)$-fold tensor product given by $\eta \mapsto 1_{C(X)} \otimes \eta$ becomes the identity map once we identify both these C$^*$-correspondences with the trivial Hilbert module $C(X)$ equipped with the left action $L_{\rho^n}$ (resp.\  $L_{\rho^{n+1}}$).  
    
    Let $\widetilde\alpha:\T(E) \longrightarrow\T(E\times_\rho X)$ be the $^*$-homomorphism constructed in Lemma~\ref{lem:induced-graph-homomorphisms}. It then follows that
    \begin{align}
       \widetilde\alpha(T_E^0(f))&=  T_{E}^0(f) \otimes 1_{C(X)},\\ 
       \widetilde\alpha(T_E^1(\xi))&= T_{E}^1(\xi) \otimes 1_{C(X)},
    \end{align}
    for all $f\in C_0(E^0)$ and $\xi\in C_d(E^1)$.
    Since $\T(E)$ is generated by the image of $T_E^0$ and $T_E^1$, we have $\widetilde\alpha(a) = a \otimes 1_{C(X)}$ for all $a \in \T(E)$.

    Recall that $I_E = \Ker(q_E)$. As shown in \cite[Section~3]{katsura2004class}, we can realise $I_E$ as the compact operators on a submodule of $C_d(E^*)$.  
    For $n\in\mathbb N_0$, write $E^n_{rg} \subseteq E^n$ for the paths of length $n$ with domain in $E^0_{rg}$. Since the domain map is continuous on paths, $E^n_{rg}$ is open in $E^n$, so we may view $C_0(E^n_{rg}) \subseteq C_0(E^n)$ and $C_d(E^n_{rg}) \subseteq C_d(E^n)$.
    We write $E^*_{rg}=\bigsqcup_{n=0}^{\infty}E^n_{rg} \subseteq E^*$ and view $C_d(E^*_{rg}) \subseteq C_d(E^*)$.
    By \cite[Proposition~3.4]{katsura2004class}, $\K(C_d(E^*_{rg})) \subseteq  \T(E)$.
    By \cite[Corollary 4.7]{katsura2004class}, we have that $I_E = \K(C_d(E^*_{rg}))$.
    
    The analogous construction also applies to $E \times_\rho X$, and therefore $I_{E \times_\rho X} = \K(C_d({(E \times_\rho X)}^*_{rg}))$. Moreover, it follows from Lemma~\ref{lem:product-regular-vertices} that the isomorphism  $C_{d}((E \times_\rho X)^*) \cong C_{d}(E^*) \boxtimes C(X)$ restricts to an isomorphism $C_{d}((E \times_\rho X)^*_{rg}) \cong C_{d}(E^*_{rg}) \boxtimes C(X)$. Hence, we have $\K(C_d({(E \times_\rho X)}^*_{rg})) \cong \K(C_{d}(E^*_{rg}) \boxtimes C(X)) \cong I_E \otimes C(X)$; see for example \cite[p.\ 37]{lance1995hilbert}.

    We summarise the conclusion of this and the previous subsection as the following proposition.
    \begin{prop}\label{prop:morph-exact-sequence}
     Let $E$ be a topological graph. 
     Let $\rho:X \rightarrow X$ be a homeomorphism of a compact Hausdorff space.
     Then there is a commuting diagram with short exact rows
         \begin{equation}\label{eqn:comparison-Toeplitz-extensions}
    \begin{tikzcd}
        0 \arrow[r]
        & I_E \arrow[r]
              \arrow[d, "\beta"']
        & \mathcal T(E) \arrow[r, "q_E"]
              \arrow[d, "\widetilde\alpha"]
        & \mathcal O(E) \arrow[r]
              \arrow[d, "\alpha"]
        & 0\\
        0 \arrow[r]
        & I_{E \times_\rho X}\arrow[r]
        & \mathcal T(E \times_\rho X) \arrow[r, "q_{E \times_\rho X}"]
        & \mathcal O(E \times_\rho X) \arrow[r]
        & 0,
    \end{tikzcd}
    \end{equation}
    where $I_E = \K(C_d(E^*_{rg}))$, $I_{E \times_\rho X} = \K(C_d((E\times_\rho X)^*_{rg})) \cong I_E \otimes C(X)$, and $\beta$ is given by $\beta(a) = a \otimes 1_{C(X)}$ for all $a \in I_E$.
    \end{prop}
 
    \subsection{The \safeK-theory computation} \label{subsec:k-theory-comparison}
    We now take $X=Z$, where $(Z,\rho)$ is the Deeley--Putnam--Strung dynamical system from Theorem~\ref{thm:DPS-Theorem}. Recall that the inclusion $\iota:\mathbb C\longrightarrow C(Z)$ given by $\lambda\mapsto\lambda 1_{C(Z)}$ is a $KK$-equivalence.

    \begin{prop}\label{prop:product-k-theory}
        Let $(Z,\rho)$ be a Deeley--Putnam--Strung minimal dynamical system as in Theorem~\ref{thm:DPS-Theorem}.
        Let $E$ be a second-countable topological graph and set $\Gamma=E\times_\rho Z$.
        The homomorphism $\alpha: \mathcal{O}(E)\rightarrow \mathcal{O}(\Gamma)$
        from Lemma~\ref{lem:induced-graph-homomorphisms} induces isomorphisms
        \begin{equation}
            \alpha_*: K_i(\mathcal O(E))
        \longrightarrow K_i(\mathcal O(\Gamma))
        \end{equation}
        for $i \in \{0,1\}$. If $E^0$ is compact, then $\alpha_*([1_{\mathcal O(E)}]_0) = [1_{\mathcal O(\Gamma)}]_0.$
    \end{prop}

\begin{proof}
    We use the commuting diagram of Toeplitz extensions
    \eqref{eqn:comparison-Toeplitz-extensions}. By naturality, there is an induced morphism of  six-term exact sequences in $K$-theory. 

    \medskip

       \begin{adjustbox}{max width=\textwidth} 
        \begin{tikzcd}
	& {K_0(I_\Gamma)} && {K_0(\mathcal T(\Gamma))} && {K_0(\mathcal{O}(\Gamma))} \\
	{K_0(I_{E})} && {K_0(\mathcal T(E))} && {K_0(\mathcal{O}(E))} \\
	& {K_1(\mathcal{O}(\Gamma))} && {K_1(\mathcal T(\Gamma))} && {K_1(I_\Gamma)} \\
	{K_1(\mathcal{O}(E))} && {K_1(\mathcal T(E))} && {K_1(I_{E})}
	\arrow[from=1-2, to=1-4]
	\arrow["{(q_\Gamma)_*}", from=1-4, to=1-6]
	\arrow[from=1-6, to=3-6]
	\arrow["{\beta_*}"{description}, from=2-1, to=1-2]
	\arrow[from=2-1, to=2-3]
	\arrow["{\tilde \alpha_*}"{description}, from=2-3, to=1-4]
	\arrow["{(q_{E})_*}", from=2-3, to=2-5]
	\arrow["{\alpha_*}"{description}, from=2-5, to=1-6]
	\arrow[from=2-5, to=4-5]
	\arrow[from=3-2, to=1-2]
	\arrow["{(q_\Gamma)_*}"', from=3-4, to=3-2]
	\arrow[from=3-6, to=3-4]
	\arrow[from=4-1, to=2-1]
	\arrow["{\alpha_*}"{description}, from=4-1, to=3-2]
	\arrow["{\tilde \alpha_*}"{description}, from=4-3, to=3-4]
	\arrow["{(q_{E})_*}"', from=4-3, to=4-1]
	\arrow["{\beta_*}"{description}, from=4-5, to=3-6]
	\arrow[from=4-5, to=4-3]
    \end{tikzcd}
    \end{adjustbox}

    \medskip

    By Lemma~\ref{lem:induced-graph-homomorphisms}, we have a commuting diagram
    \begin{equation}
        \begin{tikzcd}[column sep=large]
            C_0(E^0)
            \arrow[r, "T_E^0"]
            \arrow[d,  " - \otimes 1_{C(Z)}"]
            &
            \mathcal T(E)
            \arrow[d, "\widetilde\alpha"]
            \\
            C_0(\Gamma^0)
            \arrow[r, "T_\Gamma^0"]
            &
            \mathcal T(\Gamma),
        \end{tikzcd}
    \end{equation}
    where $T_E^0$ and $T_\Gamma^0$ come from the universal Toeplitz pairs
    and, under the identification $C_0(\Gamma^0)\cong C_0(E^0)\otimes C(Z)$, the left-hand vertical map is $f \mapsto f\otimes1_{C(Z)}$ (see the proof of Lemma~\ref{lem:induced-graph-homomorphisms}). 
    By \cite[Lemma~6.5]{katsura2004class}, both $T_E^0$ and $T_\Gamma^0$ are $KK$-equivalences. 
    By Theorem~\ref{thm:DPS-Theorem}, $f \mapsto f\otimes1_{C(Z)}$ implements an isomorphism on $K$-theory, so
    \begin{equation}
        \widetilde\alpha_*:
        K_i(\mathcal T(E))
        \longrightarrow K_i(\mathcal T(\Gamma))
    \end{equation}
    is an isomorphism for $i\in\{0,1\}$.

    Under the identification $I_\Gamma\cong I_E\otimes C(Z)$ from Proposition~\ref{prop:morph-exact-sequence}, the map $\beta$ is given by $\beta(a)=a\otimes1_{C(Z)}$ for all $a \in I_E$.
    It follows now from Theorem~\ref{thm:DPS-Theorem} that 
    \begin{equation}
        \beta_*: K_i(I_E)\longrightarrow K_i(I_\Gamma)
    \end{equation}
    is an isomorphism for $i\in\{0,1\}$.

    Applying the Five Lemma (see for example \cite[p.\ 129]{hatcher2002topology}) to the morphism of six-term exact sequences above yields that
    \begin{equation}
        \alpha_*: K_i(\mathcal O(E))
        \longrightarrow K_i(\mathcal O(\Gamma))
    \end{equation}
    is an isomorphism for $i\in\{0,1\}$.
    Finally, if $E^0$ is compact, then $\alpha$ is unital by Lemma~\ref{lem:induced-graph-homomorphisms}. 
    Consequently, $\alpha_*([1_{\mathcal O(E)}]_0) = [1_{\mathcal O(\Gamma)}]_0. $
    \end{proof}

    \section{Applying liberation to Katsura's models}\label{sec:principal-models}

    In this section, we first recall Katsura's construction of topological graph models for all unital UCT Kirchberg algebras (\cite{katsura2008class}). Then we will show that applying the technique of liberation introduced in Section~\ref{sec:dynamical-modification} still yields a topological graph algebra which is UCT Kirchberg.
    
    \subsection{Katsura's unital models} \label{subsec:construction-katsura}

    Let $G_0$ and $G_1$ be countable abelian groups and let $g\in G_0$. In \cite{katsura2008class}, Katsura constructs a topological graph $\tilde F$ such that $\mathcal O(\tilde F)$ is a unital UCT Kirchberg algebra with
    \begin{equation}
        (K_0(\mathcal O(\tilde F)), [1_{\mathcal O(\tilde F)}]_0, K_1(\mathcal O(\tilde F))) \cong (G_0,g,G_1).
    \end{equation}

    Write $M_\infty(\mathbb Z)$ for the set of column-finite integer-valued matrices indexed by $\mathbb Z_+$, and $M_\infty(\mathbb N_0)$ for its subset consisting of matrices with non-negative entries.
    Choose matrices $N\in M_\infty(\mathbb N_0)$ and $ M\in M_\infty(\mathbb Z)$
    satisfying the conditions of \cite[Lemma~4.11]{katsura2008class} for $(G_0,g,G_1)$.
    In particular, $N_{k,l}=0$ implies $M_{k,l}=0$, and the following properties hold:
    \begin{enumerate}[(i)]
        \item $N_{k,k}\geq2$ and $M_{k,k}=1$ for all $k\in\mathbb Z_+$;
        \item for all $k,l\in\mathbb Z_+$, there are $m \in \N$ and pairs $(k_i,l_i)$, $0\leq i\leq m$, with $N_{k_i,l_i}\geq1$ such that $k_0=k$, $l_i=k_{i+1}$ for $0\leq i<m$, and $l_m=l$.
    \end{enumerate}

    Define a discrete graph $E=(E^0,E^1,d_E,r_E)$ by $E^0=\mathbb Z_+$ and $ E^1=\{(k,l)\in\mathbb Z_+\times\mathbb Z_+: N_{k,l}\geq1\}$ with $d_E(k,l)=k $ and $ r_E(k,l)=l.$
    Note that property (ii) ensures that there is a path from every vertex of $E$ to every other vertex.

    Next, define $F=(F^0,F^1,d_F,r_F)$ by $F^0=E^0\times\mathbb T$ and $F^1=E^1\times\mathbb T$ with $d_F((k,l),z)=(k,z^{N_{k,l}}) $ and $r_F((k,l),z)=(l,z^{M_{k,l}}).$

    Finally, fix the vertex $w=(1,1)\in F^0$. Let $\tilde F^0=F^0\cup\{\infty\}$ be the one-point compactification of $F^0$, and set $\tilde F^1=F^1\sqcup\tilde F^0.$
    The domain and range maps agree with those of $F$ on $F^1$, while on the added copy of $\tilde F^0$ they are given by $d_{\tilde F}(v)=v $ and $
    r_{\tilde F}(v)=w$. Thus $\tilde F$ is obtained by adjoining the vertex $\infty$ and an edge from every vertex, including $\infty$, to $w$. Furthermore, note that by construction, the graph $\tilde F$ is second-countable, and its vertex space $\tilde F^0$ is compact.

    \begin{thm}[{\cite[Proposition~4.12]{katsura2008class}}] \label{thm:katsura}
        Let $\tilde F$ be the topological graph constructed above from the data $(G_0,g,G_1)$. Then $\mathcal O(\tilde F)$ is a unital UCT Kirchberg algebra with
        \begin{equation}
            (K_0(\mathcal O(\tilde F)), [1_{\mathcal O(\tilde F)}]_0, K_1(\mathcal O(\tilde F))) \cong (G_0,g,G_1).
        \end{equation}
    \end{thm}
Note that the unital Kirchberg--Phillips theorem (\cite{Ki95,Ph00}) implies that every unital UCT Kirchberg algebra is isomorphic to an algebra of this form.

    \subsection{Minimality} \label{subsec:modified-model-minimality}

    We show that $\tilde F\times_\rho X$ is minimal whenever $\tilde F$ is one of the topological graphs constructed in Subsection~\ref{subsec:construction-katsura} and $\rho$ is a minimal homeomorphism of a compact Hausdorff space $X$. The argument is an adaptation of Katsura's proof of \cite[Proposition~4.5]{katsura2008class}.

    \begin{lem}\label{lem:minimality}
        Let $E$, $F$ and $\tilde F$ be the topological graphs constructed in Subsection~\ref{subsec:construction-katsura} for some choice of $(G_0,g,G_1)$.
        Let $a\in\tilde F^0$, let $V\subseteq\tilde F^0$ be a non-empty open set, and let $S\subseteq\mathbb N$ be infinite. Then there are $n\in S$ and $\mu\in\tilde F^n$ such that
        \begin{equation}
            d_{\tilde F}(\mu)=a \quad \text{and} \quad  r_{\tilde F}(\mu)\in V.
        \end{equation}
    \end{lem}

    \begin{proof}
        Recall that $F^0=\mathbb Z_+\times\mathbb T$ is dense in $\tilde F^0$. Thus there are $l\in\mathbb Z_+$ and a non-empty open arc $B\subseteq\mathbb T$ such that
        \begin{equation}
            \{l\}\times B\subseteq V.
        \end{equation}
        
        First suppose that $a=(k,t)\in F^0$. By property (ii) in Subsection~\ref{subsec:construction-katsura}, there is a path $\nu\in E^m$, for some $m\geq1$, with $d_E(\nu)=k$ and $r_E(\nu)=l$.
        We can obtain a lift $\nu^\prime\in F^m$ with $d_F(\nu^\prime)=(k,t)$ and $ r_F(\nu^\prime)=(l,s)$ for some $s\in\mathbb T$ by choosing the circle coordinates successively, starting at the domain of $\nu$. For this, we use that each map $z\mapsto z^{N_{k',l'}}$ associated to an edge $(k',l')$ of $E$ is surjective.
        
        Next, set $q=N_{l,l}\geq2$ and recall that $M_{l,l}=1$. Thus, for every $j\geq1$ and every $u\in\mathbb T$ with $u^{q^j}=s$, there exists a path
        \begin{equation}
            \lambda_{j,u} = \bigl( ((l,l),u), ((l,l),u^q), \ldots, ((l,l),u^{q^{j-1}})
            \bigr)
        \end{equation}
        which has domain $(l,s)$ and range $(l,u)$. 
        The $q^j$-many solutions of $u^{q^j}=s$ are equally spaced around the circle.
        Hence, for every sufficiently large $j$, one of these solutions lies in $B$. Choose such a solution $u$. Then the concatenated path
        \begin{equation}
            (\lambda_{j,u},\nu^\prime)\in F^{m+j}
        \end{equation}
        exists and has domain $a=(k,t)$ and range in $\{l\}\times B\subseteq V$. Since $S$ is infinite, we can choose $n=m+j\in S$ with $j$ sufficiently large. This proves the claim for $a\in F^0$.

        Now suppose that $a=\infty$. Recall that there is an added edge $e_\infty\in\tilde F^1$ with $d_{\tilde F}(e_\infty)=\infty $ and $ r_{\tilde F}(e_\infty)=w$
        and that the set $S'=\{n-1:n\in S,\, n\geq2\}$ is still infinite.
        Applying the first part to $w$, $V$ and $S'$ gives $m\in S'$ and a path $\eta\in\tilde F^m$ with domain $w$ and range in $V$.
        Then $(\eta,e_\infty)$ has length $m+1\in S$, domain $\infty$ and range in $V$.
    \end{proof}

    \begin{prop}\label{prop:minimality}
        Let $\tilde F$ be the topological graph constructed in Subsection~\ref{subsec:construction-katsura} for some choice of $(G_0,g,G_1)$ and let $\rho: X\to X$ be a minimal homeomorphism of a compact Hausdorff space $X$. Then $\tilde F\times_\rho X$ is minimal.
    \end{prop}

    \begin{proof}
        Set $\Gamma=\tilde F\times_\rho X$. It suffices to show that every positive orbit is dense. Fix $(v,x)\in\Gamma^0$ and a non-empty basic open set $U\times V\subseteq\tilde F^0\times X$.

        Since $\rho$ is minimal and $X$ is compact, every forward orbit is dense, and the set
        \begin{equation}
            S=\{n\in\mathbb N:\rho^n(x)\in V\}
        \end{equation}
        is infinite.
        By Lemma~\ref{lem:minimality}, there are $n\in S$ and $\mu\in\tilde F^n$ such that $d_{\tilde F}(\mu)=v$ and $r_{\tilde F}(\mu)\in U.$
        Using the path space homeomorphism from Lemma~\ref{lem:path-homeo}, set $\lambda=\Phi_n(\mu,\rho^n(x)).$
        Then, by equation~\eqref{eqn:modified-path-domain-range}, $d_\Gamma(\lambda)=(v,x)$ and $r_\Gamma(\lambda) = (r_{\tilde F}(\mu),\rho^n(x)) \in U\times V.$
        Thus every positive orbit is dense, and $\Gamma$ is minimal.
    \end{proof}

    \subsection{Contraction} \label{subsec:modified-model-contraction}

    Next, we show that $\tilde F\times_\rho X$ is contracting whenever $\tilde F$ is one of the topological graphs constructed in Subsection~\ref{subsec:construction-katsura} and $\rho$ is a minimal homeomorphism of a compact Hausdorff space $X$. 
    
    \begin{prop}\label{prop:contracting}
        Let $\tilde F$ be the topological graph constructed in Subsection~\ref{subsec:construction-katsura} for some choice of $(G_0,g,G_1)$ and let $\rho: X\to X$ be a minimal homeomorphism of a compact Hausdorff space $X$. Then $\tilde F\times_\rho X$ is contracting.
    \end{prop}

    \begin{proof}
        Set $\Gamma=\tilde F\times_\rho X$ and fix $x_0\in X$. We show that $\Gamma$ is contracting at $(w,x_0)$, where $w=(1,1)\in F^0$.
        First note that, by the proof of Proposition~\ref{prop:minimality},
        \begin{equation}
            \overline{\operatorname{Orb}^+((w,x_0))} = \Gamma^0.
        \end{equation}

        Now, let $W$ be a neighbourhood of $(w,x_0)$. Choose an open neighbourhood $U$ of $x_0$ in $X$ and an open arc $J\subseteq\mathbb T$ containing $1$, with $\overline J\neq\mathbb T$, such that
        \begin{equation}
            V:=(\{1\}\times J)\times U\subseteq W.
        \end{equation}
        We show that $V$ is contracting. First note that its closure is compact because $X$ is compact.

        Set $q=N_{1,1}\geq2$, and let $|J|$ denote the length of the arc $J$. Choose $L\geq1$ big enough such that
        \begin{equation}
            4\pi\sum_{n=L}^{\infty}q^{-n}<|J|.
        \end{equation}
        We can then choose pairwise disjoint open arcs $J_n\subseteq J$, for $n\geq L$, such that
        \begin{equation}
            \frac{2\pi}{q^n}<|J_n|<\frac{4\pi}{q^n}.
        \end{equation}
        Since $|J_n|>2\pi/q^n$, the map $t\mapsto t^{q^n}$ sends $J_n$ onto $\mathbb T$.

        By minimality of $\rho$ and compactness of $X$, we can find $L'\geq L$ with
        \begin{equation}
            \bigcup_{n=L}^{L'}\rho^{-n}(U)=X.
        \end{equation}

        For $n\geq1$ and $t\in\mathbb T$, define
        \begin{equation}
        \lambda_n(t) = \bigl( ((1,1),t), ((1,1),t^q), \ldots, ((1,1),t^{q^{n-1}}) \bigr) \in F^n \subset \tilde{F}^n.
        \end{equation}
        This is a path because $M_{1,1}=1$ and $N_{1,1} = q.$ For $L\leq n\leq L'$, we then set
        \begin{equation}
            U_n = \{\Phi_n(\lambda_n(t),x): t\in J_n,\, x\in U\} \subseteq\Gamma^n,
        \end{equation}
        using the path space homeomorphisms from Lemma~\ref{lem:path-homeo}.
        The set of paths in $\tilde F^n$ whose edges all lie in $\{(1,1)\}\times\mathbb T$ is open and is parametrised (homeomorphically) exactly by $t\mapsto\lambda_n(t)$.
        Hence $\lambda_n(J_n)$ is open in $\tilde F^n$, and since $\Phi_n$ is a homeomorphism, $U_n$ is a non-empty open subset of $\Gamma^n$.

        By equation~\eqref{eqn:modified-path-domain-range},
        \begin{equation}\begin{split}
            r_\Gamma(\Phi_n(\lambda_n(t),x)) &=((1,t),x),\\
            d_\Gamma(\Phi_n(\lambda_n(t),x)) &=((1,t^{q^n}),\rho^{-n}(x)).
        \end{split}\end{equation}
        Consequently, $r_\Gamma(U_n) = (\{1\}\times J_n)\times U \subseteq V$ and $d_\Gamma(U_n) = (\{1\}\times\mathbb T)\times\rho^{-n}(U).$
    
        Moreover, $U_n\pitchfork U_m=\emptyset$ whenever $n\neq m$, since the circle coordinates of the first edges belong to the disjoint arcs $J_n$ and $J_m$.
        Finally,
        \begin{equation}
            \overline V = (\{1\}\times\overline J)\times\overline U \subsetneq (\{1\}\times\mathbb T)\times X = \bigcup_{n=L}^{L'}d_\Gamma(U_n).
        \end{equation}
        Thus $V$ is contracting, and $\Gamma$ is contracting at $(w,x_0)$.
    \end{proof}

\subsection{Dimension of the boundary path space}
\label{subsec:boundary-path-dimension}
    We can bound the dimension of the unit space of the boundary path groupoid associated to $\tilde F\times_\rho X$ whenever $\tilde F$ is one of the topological graphs constructed in Subsection~\ref{subsec:construction-katsura} and $\rho:X \rightarrow X$ is a homeomorphism of a compact metric space of finite covering dimension.
    
    \begin{prop}\label{prop:spectrum-dimension}
         Let $\tilde F$ be the topological graph constructed in Subsection~\ref{subsec:construction-katsura} for some choice of $(G_0,g,G_1)$.
        Let $X$ be a compact metric space, and let $\rho: X\to X$ be a homeomorphism.
        Let $\G$ denote the boundary path groupoid associated to $\tilde F\times_\rho X$ (see Subsection~\ref{subsec:top-graphs-groupoids}).
        Then
        \begin{equation}
            \dim(\G^{(0)})\leq2\dim(X)+3.
        \end{equation}
        In particular, if $X$ has finite covering dimension, then $\G^{(0)}$ also has finite covering dimension.
    \end{prop}

    \begin{proof}
        Set $\Gamma=\tilde F\times_\rho X$. Since $\Gamma$ is second-countable, $\partial\Gamma$ is second-countable, locally compact, Hausdorff and metrizable.
        Its topology is therefore determined by the convergence criterion in Lemma~\ref{lem:Kumjian-Li}.

        Recall that taking the one-point compactification of a locally compact separable metrizable space does not change its covering dimension. Thus, $\dim(\tilde F^0)=1$.
        By \cite[Theorem~3.2.14]{engelking1978dimension},
        \begin{equation}\label{eqn:modified-vertex-dimension}
            \dim(\Gamma^0) = \dim(\tilde F^0\times X) \leq \dim(\tilde F^0)+\dim(X) = 1+\dim(X).
        \end{equation}
        The product theorem applies since both factors are compact metrizable spaces and hence normal and strongly paracompact.

        We first estimate the dimensions of the finite path spaces.
        For each $n\in\mathbb N_0$, the domain map $d_\Gamma: \Gamma^n\to\Gamma^0$ is a local homeomorphism.
        Since $\Gamma^n$ is second-countable, it has a countable open cover $(U_i)_{i=1}^{\infty}$ such that each restriction $d_\Gamma|_{U_i}$ is a homeomorphism onto an open subset of $\Gamma^0$.

        As $\Gamma^n$ is metrizable, every $U_i$ is a countable union of sets which are closed in $\Gamma^n$. Each of these closed sets is homeomorphic to a closed subspace of $\Gamma^0$ and therefore has dimension at most $\dim(\Gamma^0)$ by \cite[Theorem~3.1.19]{engelking1978dimension}.
        Here the subspace theorem applies because $\Gamma^0$ is separable and metrizable.
        Consequently, $\Gamma^n$ is a countable union of closed sets of dimension at most $\dim(\Gamma^0)$. The countable closed sum theorem (\cite[Theorem~3.1.8]{engelking1978dimension}), applicable since $\Gamma^n$ is metrizable, gives $\dim(\Gamma^n)\leq\dim(\Gamma^0).$

        Recall that $\partial\Gamma = \Gamma^\infty \sqcup \{\mu\in\Gamma^*:  d_\Gamma(\mu) \in\Gamma^0_{sg}\}$ as a set and has topology defined in Section \ref{subsec:graphs-and-paths}. We write $B=
        \{\mu\in\Gamma^*: d_\Gamma(\mu) \in \Gamma^0_{sg}\} $ and note that $B$ carries the subspace topology inherited from $\partial\Gamma$, rather than the disjoint union topology on $\Gamma^*$.
        Furthermore, for $n\in\mathbb N_0$, we set $ B_n= \{\mu\in\Gamma^n: d_\Gamma(\mu) \in\Gamma^0_{sg}\}.$
        
        Since $\Gamma^0_{sg}$ is closed, $B_n$ is a closed subset of $\Gamma^n$.
        Moreover, Lemma~\ref{lem:Kumjian-Li} shows that its topology as a subspace of $\partial\Gamma$ agrees with its topology as a subspace of $\Gamma^n$. Indeed, for paths of fixed length $n$, convergence to a path of length $n$ is exactly convergence in $\Gamma^n$.

        Since $B_n$ is closed in $\Gamma^n$,  \cite[Theorem~3.1.19]{engelking1978dimension} gives $\dim(B_n)\leq\dim(\Gamma^n) \leq\dim(\Gamma^0) \le 1 + \dim(X)$ for all $n \in \N_0$.
        Each $B_n$ is locally compact and second-countable, so we can write it as a countable union of compact sets. 
        Since $\partial\Gamma$ is Hausdorff, these compact sets are also closed in both $B_n$ and $B$, and each has dimension at most $\dim(B_n) = 1 + \dim(X)$. 
        Thus, $B=\bigcup_{n=0}^\infty B_n$ is a countable union of closed subsets of dimension at most $1 + \dim(X)$. The countable closed sum theorem (\cite[Theorem~3.1.8]{engelking1978dimension}) applies because $B$ is metrizable, and gives $\dim(B)\leq1 + \dim(X).$
    
        Next, we consider the infinite path space $\Gamma^\infty$. For $n\geq1$, let $p_n:\Gamma^\infty \to\Gamma^n$ be given by $p_n((e_i)_{i=1}^{\infty})=(e_1,\ldots,e_n)$ and set $P_n=p_n(\Gamma^\infty).$
        Thus, $P_n$ is the space of paths of length $n$ which are extendable to infinite paths. We define surjections $q_n: P_{n+1}\to  P_n$ by $q_n(e_1,\ldots,e_{n+1})=(e_1,\ldots,e_n)$ and note that they are continuous.
        There is a canonical bijection
        \begin{equation}\begin{split}
            \Gamma^\infty &\longrightarrow\varprojlim(P_n,q_n),\\
            (e_1,e_2,\ldots) &\longmapsto ((e_1),(e_1,e_2),(e_1,e_2,e_3),\ldots),
        \end{split}\end{equation}
        which is a homeomorphism.
        Indeed, by Lemma~\ref{lem:Kumjian-Li}, a sequence of infinite paths converges to an infinite path if and only if all its finite initial segments converge.
        This is precisely convergence in the inverse-limit topology.
        In particular, both spaces are metrizable.

        Each $P_n$ is a separable metrizable subspace of $\Gamma^n$. Hence \cite[Theorem~3.1.19]{engelking1978dimension} gives
        \begin{equation}
            \dim(P_n)\leq\dim(\Gamma^n) \leq\dim(\Gamma^0)\leq 1 + \dim(X).
        \end{equation}
        Applying the dimension theorem for countable inverse limits of separable metrizable spaces (see \cite[Theorem 2]{Nagami59} or \cite[Exercise~7.3.I]{engelking1989General}), we obtain
        \begin{equation}\label{eqn:infinite-path-dimension}
            \dim(\Gamma^\infty)\leq\dim(\Gamma^0) \leq 1 + \dim(X).
        \end{equation}

        Since $\partial\Gamma$ is metrizable, the union inequality \cite[Theorem~3.1.17]{engelking1978dimension} gives
        \begin{equation} 
            \dim(\partial\Gamma) \leq \dim(B)+\dim(\Gamma^\infty)+1.
        \end{equation}
        By the previous computations, this yields
        \begin{equation}
            \dim(\G^{(0)}) =\dim(\partial\Gamma) \leq 2\dim(X)+3. \qedhere
        \end{equation}
    \end{proof}

\section{Classifiable \cstar-algebras, groupoid models and \cstar-diagonals}\label{sec:main-proofs}
    We now have all the ingredients needed to prove the main theorems discussed in the introduction. For the benefit of the reader, we restate these results here too.

    \begin{thm} [Theorem~\ref{Thm:Kirchberg}]
        Let $A$ be a Kirchberg algebra satisfying the UCT. Then there exists a principal, second-countable, locally compact, Hausdorff, étale groupoid $\G$ with a finite-dimensional unit space such that $A \cong C_r^*(\G)$.
    \end{thm}
    
    \begin{proof}
        By \cite[Theorem~1.2]{Zhang92}, a Kirchberg algebra is either unital or stable.
        If $A$ is stable, the result was proven in \cite[Theorem~A]{ES26Stable}. 
        
        Let $A$ be a unital UCT Kirchberg algebra. Set $G_0=K_0(A), G_1 =K_1(A)$ and $g=[1_A]_0$.
        Choose the graph $\tilde F$ from Theorem~\ref{thm:katsura} for this $K$-theory.
        Let $(Z,\rho)$ be the Deeley--Putnam--Strung dynamical system from Theorem~\ref{thm:DPS-Theorem}, and set $\Gamma=\tilde F\times_\rho Z.$ 
        The graph $\Gamma$ is second-countable, minimal by Proposition~\ref{prop:minimality}, and contracting by Proposition~\ref{prop:contracting}.
        Thus $\mathcal O(\Gamma)$ is a UCT Kirchberg algebra by \cite[Corollary~B]{katsura2008class}. Moreover, $\Gamma^0=\tilde F^0\times Z$ is compact, so $\mathcal O(\Gamma)$ is unital.

        By Proposition~\ref{prop:product-k-theory}, the algebras $\mathcal O(\tilde F)$ and $\mathcal O(\Gamma)$ have isomorphic $K$-theory, with the isomorphism on $K_0$ preserving the class of the unit.
        The unital Kirchberg--Phillips theorem (\cite{Ki95,Ph00}) therefore gives $\mathcal O(\Gamma) \cong \mathcal O(\tilde F) \cong A.$

        By Proposition~\ref{prop:liberation-DPS}, the boundary path groupoid $\G=\G(\partial\Gamma,\sigma_\Gamma)$  is principal, second-countable, locally compact, Hausdorff and étale, and $\mathcal O(\Gamma)\cong C_r^*(\G)$.
        By Proposition~\ref{prop:spectrum-dimension}, $\G^{(0)}=\partial\Gamma$ has finite covering dimension.
    \end{proof}

    \begin{thm} [Theorem~\ref{Thm:every-classifiable-diagonal}]
        Every simple, separable, nuclear, UCT, Toms--Winter regular C$^*$-algebra $A$ has a C$^*$-diagonal $D \subseteq A$. Moreover, $D$ can be chosen to have finite-dimensional spectrum.
    \end{thm}
    \begin{proof}
        Let $A$ be a simple, separable, nuclear, UCT, Toms--Winter regular C$^*$-algebra. If $A$ is elementary, there is a standard C$^*$-diagonal with zero-dimensional spectrum: either $\C^n \subseteq M_n$ or $c_0 \subseteq \K$. Hence, we may also assume that $A$ is non-elementary.

        Since $A$ is Toms--Winter regular, it is Jiang--Su stable. Hence, $A$ is either stably finite or purely infinite (see for example \cite[Theorem 4.1.10(ii)]{rordam2002classification}).
        
        If $A$ is purely infinite, it is a UCT Kirchberg algebra. By Theorem~\ref{Thm:Kirchberg}, it has a principal groupoid model with finite-dimensional unit space. Kumjian--Renault theory therefore gives a C$^*$-diagonal in $A$ with finite-dimensional spectrum.
        
        If $A$ is stably finite, the existence of a C$^*$-diagonal follows from Xin Li's construction \cite[Theorems 1.2 and 1.3]{li2020every} and the classification theorems of the Elliott classification programme (\cite{EGLN25, GLN20a, GLN20b, TWW17, GL23, CGSTW}), noting that the required classification results in the non-unital stably projectionless case now exist (see for example \cite[Theorem~14.39]{gong2024classification}). By \cite[Corollary~1.8]{li2020every}, the C$^*$-diagonal can be chosen to have finite-dimensional spectrum. 
    \end{proof}

    \begin{cor}[Corollary~\ref{Cor:UCT-diagonal}]
        Let $A$ be a simple, separable, nuclear, Toms--Winter regular C$^*$-algebra. Then the following are equivalent:
        \begin{enumerate}[(i)]
            \item $A$ satisfies the UCT;
            \item $A$ has a C$^*$-diagonal;
            \item $A$ has a Cartan subalgebra.
        \end{enumerate}
    \end{cor}
    \begin{proof}
        By Theorem~\ref{Thm:every-classifiable-diagonal}, we have (i)$\Rightarrow$(ii). 
        A C$^*$-diagonal is a Cartan subalgebra by definition, so (ii)$\Rightarrow$(iii). By \cite[Corollary~1.2]{barlak2017cartan}, we have (iii)$\Rightarrow$(i).
    \end{proof}


\begin{thebibliography}{10}

\bibitem{barlak2017cartan}
S.~Barlak and X.~Li.
\newblock Cartan subalgebras and the {UCT} problem.
\newblock {\em Adv. Math.}, 316:748--769, 2017.

\bibitem{brown2016purely}
J.~Brown, L.~Clark, A.~Sierakowski, and A.~Sims.
\newblock Purely infinite simple \textnormal{C}$^*$-algebras that are principal
  groupoid \textnormal{C}$^*$-algebras.
\newblock {\em J. Math. Anal. Appl.}, 439(1):213--234, 2016.

\bibitem{CGSTW}
J.~Carri\'on, J.~Gabe, C.~Schafhauser, A.~Tikuisis, and S.~White.
\newblock Classifying {$^*$}-homomorphisms {I}: Unital, simple, nuclear
  {$\rm{C}^*$}-algebras.
\newblock {\em Acta Math.}
\newblock to appear, arXiv:2307.06480, 2026.

\bibitem{CFH20}
L.~Clark, J.~Fletcher, and A.~an~Huef.
\newblock All classifiable {K}irchberg algebras are \textnormal{C}$^*$-algebras
  of ample groupoids.
\newblock {\em Expo. Math.}, 38(4):559--565, 2020.

\bibitem{cuntz1977simple}
J.~Cuntz.
\newblock Simple \textnormal{C}$^*$-algebras generated by isometries.
\newblock {\em Comm. Math. Phys.}, 57(2):173--185, 1977.

\bibitem{deaconu1995groupoids}
V.~Deaconu.
\newblock Groupoids associated with endomorphisms.
\newblock {\em Trans. Amer. Math. Soc.}, 347(5):1779--1786, 1995.

\bibitem{deeley2018constructing}
R.~Deeley, I.~Putnam, and K.~Strung.
\newblock Constructing minimal homeomorphisms on point-like spaces and a
  dynamical presentation of the {Jiang--Su} algebra.
\newblock {\em J. Reine Angew. Math.}, 742:241--261, 2018.

\bibitem{deeley2023min-homeo}
R.~Deeley, I.~Putnam, and K.~Strung.
\newblock Minimal homeomorphisms and topological {$K$}-theory.
\newblock {\em Groups Geom. Dyn.}, 17(2):501--532, 2023.

\bibitem{elliott1976classification}
G.~Elliott.
\newblock On the classification of inductive limits of sequences of semisimple
  finite-dimensional algebras.
\newblock {\em J. Algebra}, 38(1):29--44, 1976.

\bibitem{El95}
G.~Elliott.
\newblock The classification problem for amenable \textnormal{C}$^*$-algebras.
\newblock In {\em Proceedings of the {I}nternational {C}ongress of
  {M}athematicians, {V}ol.\ 1, 2 ({Z}\"urich, 1994)}, pages 922--932.
  Birkh\"auser, Basel, 1995.

\bibitem{El96}
G.~Elliott.
\newblock An invariant for simple {$\rm{C}^*$}-algebras.
\newblock In {\em Canadian {M}athematical {S}ociety. 1945--1995, {V}ol. 3},
  pages 61--90. Canadian Math. Soc., Ottawa, ON, 1996.

\bibitem{EGLN25}
G.~Elliott, G.~Gong, H.~Lin, and Z.~Niu.
\newblock On the classification of simple amenable {$\rm {C}^*$}-algebras with
  finite decomposition rank, {II}.
\newblock {\em J. Noncommut. Geom.}, 19(1):73--104, 2025.

\bibitem{elliott2016k}
G.~Elliott and A.~Sierakowski.
\newblock {$K$}-theory of certain purely infinite crossed products.
\newblock {\em J. Math. Anal. Appl.}, 443(1):409--430, 2016.

\bibitem{engelking1978dimension}
R.~Engelking.
\newblock {\em Dimension theory}, volume~19 of {\em North-Holland Mathematical
  Library}.
\newblock North-Holland Publishing Co., Amsterdam-Oxford-New York; PWN---Polish
  Scientific Publishers, Warsaw, 1978.
\newblock Translated from the Polish and revised by the author.

\bibitem{engelking1989General}
R.~Engelking.
\newblock {\em General topology}, volume~6 of {\em Sigma Series in Pure
  Mathematics}.
\newblock Heldermann Verlag, Berlin, second edition, 1989.
\newblock Translated from the Polish by the author.

\bibitem{TomsWinter}
S.~Evington, S.~Hua, C.~Schafhauser, A.~Seth, A.~Tikuisis, and S.~White.
\newblock Regular simple nuclear \textnormal{C}$^*$-algebras.
\newblock arXiv:2610.09390, 2026.

\bibitem{ES26preprint}
S.~Evington and P.~Sibbel.
\newblock Principal groupoid models for {Cuntz} algebras and their dynamic
  asymptotic dimension.
\newblock arXiv:2512.01812, 2025.

\bibitem{ES26}
S.~Evington and P.~Sibbel.
\newblock {$\rm{C}^*$}-diagonals with {C}antor spectrum in {C}untz algebras.
\newblock {\em J. Funct. Anal.}, 290(12):Paper No. 111418, 2026.

\bibitem{ES26Stable}
S.~Evington and P.~Sibbel.
\newblock Principal groupoid models for stable {UCT} {Kirchberg} algebras.
\newblock arXiv:2605.30147, 2026.

\bibitem{GGNV24}
E.~Gardella, S.~Geffen, P.~Naryshkin, and A.~Vaccaro.
\newblock Dynamical comparison and {$\mathcal{Z}$}-stability for crossed
  products of simple {C$^*$}-algebras.
\newblock {\em Adv. Math.}, 438:Paper No. 109471, 31, 2024.

\bibitem{GardellaPerera26}
E.~Gardella and F.~Perera.
\newblock The modern theory of {C}untz semigroups of {$C^{\ast}$}-algebras.
\newblock {\em EMS Surv. Math. Sci.}, 13(1):133--214, 2026.

\bibitem{GL23}
G.~Gong and H.~Lin.
\newblock A review of the {Elliott program} of classification of simple
  amenable {$\rm{C}^*$}-algebras.
\newblock arXiv:2311.14238, 2023.

\bibitem{gong2024classification}
G.~Gong and H.~Lin.
\newblock On classification of nonunital amenable simple {$\rm{C}^*$}-algebras,
  {IV}: {Stably} projectionless {$\rm{C}^*$}-algebras.
\newblock {\em Ann. K-Theory}, 9(2):143--339, 2024.

\bibitem{GLN20a}
G.~Gong, H.~Lin, and Z.~Niu.
\newblock A classification of finite simple amenable {$\mathcal{Z}$}-stable
  {$\rm{C}^*$}-algebras, {I}: {C$^*$}-algebras with generalized tracial rank
  one.
\newblock {\em C. R. Math. Acad. Sci. Soc. R. Can.}, 42(3):63--450, 2020.

\bibitem{GLN20b}
G.~Gong, H.~Lin, and Z.~Niu.
\newblock A classification of finite simple amenable {$\mathcal{Z}$}-stable
  {$\rm{C}^*$}-algebras, {II}: {C$^*$}-algebras with rational generalized
  tracial rank one.
\newblock {\em C. R. Math. Acad. Sci. Soc. R. Can.}, 42(4):451--539, 2020.

\bibitem{hatcher2002topology}
A.~Hatcher.
\newblock {\em Algebraic topology}.
\newblock Cambridge University Press, Cambridge, 2002.

\bibitem{hjelmborgstability}
J.~Hjelmborg.
\newblock {\em On Stability and Pure Infiniteness of
  \textnormal{C}$^*$-Algebras}.
\newblock PhD thesis, University of Southern Denmark, 2000.
\newblock Available at https://web.math.ku.dk/rordam/students/jh-thesis.ps.

\bibitem{JS99}
X.~Jiang and H.~Su.
\newblock On a simple unital projectionless {$\rm{C}^*$}-algebra.
\newblock {\em Amer. J. Math.}, 121(2):359--413, 1999.

\bibitem{katsura2004class}
T.~Katsura.
\newblock A class of \textnormal{C}$^*$-algebras generalizing both graph
  algebras and homeomorphism \textnormal{C}$^*$-algebras {I}, fundamental
  results.
\newblock {\em Trans. Amer. Math. Soc.}, 356(11):4287--4322, 2004.

\bibitem{katsura2006class2}
T.~Katsura.
\newblock A class of \textnormal{C}$^*$-algebras generalizing both graph
  algebras and homeomorphism \textnormal{C}$^*$-algebras {II}, examples.
\newblock {\em Internat. J. Math.}, 17(7):791--833, 2006.

\bibitem{katsura2008class}
T.~Katsura.
\newblock A class of \textnormal{C}$^*$-algebras generalizing both graph
  algebras and homeomorphism \textnormal{C}$^*$-algebras {IV}, pure
  infiniteness.
\newblock {\em J. Funct. Anal.}, 254(5):1161--1187, 2008.

\bibitem{KerrSzabo20}
D.~Kerr and G.~Szab\'o.
\newblock Almost finiteness and the small boundary property.
\newblock {\em Comm. Math. Phys.}, 374(1):1--31, 2020.

\bibitem{Ki95}
E.~Kirchberg.
\newblock Exact \textnormal{C}$^*$-algebras, tensor products, and the
  classification of purely infinite algebras.
\newblock In {\em Proceedings of the {I}nternational {C}ongress of
  {M}athematicians, {V}ol.\ 1, 2 ({Z}\"urich, 1994)}, pages 943--954.
  Birkh\"auser, Basel, 1995.

\bibitem{kumjian1986c}
A.~Kumjian.
\newblock On \textnormal{C}$^*$-diagonals.
\newblock {\em Canad. J. Math.}, 38(4):969--1008, 1986.

\bibitem{kumjian2017twisted}
A.~Kumjian and H.~Li.
\newblock Twisted topological graph algebras are twisted groupoid
  \textnormal{C}$^*$-algebras.
\newblock {\em J. Operator Theory}, 78(1):201--225, 2017.

\bibitem{lance1995hilbert}
E.~C. Lance.
\newblock {\em Hilbert \textnormal{C}$^*$-modules. A toolkit for operator
  algebraists}, volume 210 of {\em London Math. Soc. Lecture Note Ser.}
\newblock Cambridge University Press, Cambridge, 1995.

\bibitem{LLW23}
K.~Li, H.~Liao, and W.~Winter.
\newblock The diagonal dimension of sub-\textnormal{C}$^*$-algebras.
\newblock {\em Proc. Lond. Math. Soc. (3)}, 133(4):Paper No. e70234, 2026.

\bibitem{li2020every}
X.~Li.
\newblock Every classifiable simple \textnormal{C}$^*$-algebra has a {Cartan}
  subalgebra.
\newblock {\em Invent. Math.}, 219(2):653--699, 2020.

\bibitem{li2019cartan}
X.~Li and J.~Renault.
\newblock Cartan subalgebras in \textnormal{C}$^*$-algebras. {Existence} and
  uniqueness.
\newblock {\em Trans. Amer. Math. Soc.}, 372(3):1985--2010, 2019.

\bibitem{MatuiSato12}
H.~Matui and Y.~Sato.
\newblock {$\mathcal{Z}$}-stability of crossed products by strongly outer
  actions.
\newblock {\em Comm. Math. Phys.}, 314(1):193--228, 2012.

\bibitem{MatuiSato14}
H.~Matui and Y.~Sato.
\newblock {$\mathcal{Z}$}-stability of crossed products by strongly outer
  actions {II}.
\newblock {\em Amer. J. Math.}, 136(6):1441--1496, 2014.

\bibitem{Nagami59}
K.~Nagami.
\newblock Finite-to-one closed mappings and dimension. {II}.
\newblock {\em Proc. Japan Acad.}, 35:437--439, 1959.

\bibitem{OAI:Toms-Winter}
{OpenAI}.
\newblock {Tracial projection methods and uniform property $\Gamma$}.
\newblock OpenAI Math Release preprint
  \href{https://github.com/openai/math/blob/main/preprints/Tracial-projection-methods-and-uniform-property-Gamma-September-23-2026/paper.pdf}{OAI:Tracial-https://github.com/openai/math},
  2026.

\bibitem{Ph00}
N.~C. Phillips.
\newblock A classification theorem for nuclear purely infinite simple
  \textnormal{C}$^*$-algebras.
\newblock {\em Doc. Math.}, 5:49--114, 2000.

\bibitem{renault1980groupoid}
J.~Renault.
\newblock {\em A groupoid approach to \textnormal{C}$^*$-algebras}, volume 793
  of {\em Lecture Notes in Math.}
\newblock Springer, Berlin, 1980.

\bibitem{renault2000cuntz}
J.~Renault.
\newblock Cuntz-like algebras.
\newblock In {\em Operator theoretical methods ({T}imi\c soara, 1998)}, pages
  371--386. The Theta Foundation, Bucharest, 2000.

\bibitem{renault2008cartan}
J.~Renault.
\newblock {Cartan} subalgebras in \textnormal{C}$^*$-algebras.
\newblock {\em Irish Math. Soc. Bull.}, 61:29--63, 2008.

\bibitem{rieffel1981RotationAlgebras}
M.~Rieffel.
\newblock \textnormal{C}$^*$-algebras associated with irrational rotations.
\newblock {\em Pacific J. Math.}, 93(2):415--429, 1981.

\bibitem{rordam2002classification}
M.~R{\o}rdam.
\newblock Classification of nuclear, simple \textnormal{C}$^*$-algebras.
\newblock In {\em Classification of nuclear \textnormal{C}$^*$-Algebras.
  {Entropy} in operator algebras}, volume 126 of {\em Encyclopaedia Math.
  Sci.}, pages 1--145. Oper. Alg. Non-commut. Geom., 7, Springer, Berlin, 2002.

\bibitem{rordam2012purely}
M.~R{\o}rdam and A.~Sierakowski.
\newblock Purely infinite \textnormal{C}$^*$-algebras arising from crossed
  products.
\newblock {\em Ergodic Theory Dynam. Systems}, 32(1):273--293, 2012.

\bibitem{RS87}
J.~Rosenberg and C.~Schochet.
\newblock The {K}{\"u}nneth theorem and the universal coefficient theorem for
  {K}asparov's generalized {$K$}-functor.
\newblock {\em Duke Math. J.}, 55(2):431--474, 1987.

\bibitem{STW25_99problems}
C.~Schafhauser, A.~Tikuisis, and S.~White.
\newblock Nuclear \textnormal{C}$^*$-algebras: 99 problems.
\newblock {\em M\"unster J. Math.}
\newblock to appear, arXiv:2506.10902, 2026.

\bibitem{sibbel2024cantor}
P.~Sibbel and W.~Winter.
\newblock A {C}antor spectrum diagonal in $\mathcal{O}_2$.
\newblock {\em Proc. Amer. Math. Soc. Ser. B}, 12:210--217, 2025.

\bibitem{spielberg2007graph}
J.~Spielberg.
\newblock Graph-based models for {K}irchberg algebras.
\newblock {\em J. Operator Theory}, 57(2):347--374, 2007.

\bibitem{suzuki2017amenable}
Y.~Suzuki.
\newblock Amenable minimal {C}antor systems of free groups arising from
  diagonal actions.
\newblock {\em J. Reine Angew. Math.}, 722:183--214, 2017.

\bibitem{Thiel26}
H.~Thiel.
\newblock Ranks of operators and pureness of \textnormal{C}$^*$-algebras.
\newblock arXiv:2610.06761, 2026.

\bibitem{TWW17}
A.~Tikuisis, S.~White, and W.~Winter.
\newblock Quasidiagonality of nuclear {$\rm{C}^*$}-algebras.
\newblock {\em Ann. of Math. (2)}, 185(1):229--284, 2017.

\bibitem{Tu99}
J.~Tu.
\newblock La conjecture de {B}aum-{C}onnes pour les feuilletages moyennables.
\newblock {\em $K$-Theory}, 17(3):215--264, 1999.

\bibitem{WhiteICM}
S.~White.
\newblock Abstract classification theorems for amenable
  \textnormal{C}$^*$-algebras.
\newblock In {\em I{CM}---{I}nternational {C}ongress of {M}athematicians.
  {V}ol. 4. {S}ections 5--8}, pages 3314--3338. EMS Press, Berlin, 2023.

\bibitem{winter2018structure}
W.~Winter.
\newblock Structure of nuclear \textnormal{C}$^*$-algebras: from
  quasidiagonality to classification and back again.
\newblock In {\em Proceedings of the {I}nternational {C}ongress of
  {M}athematicians---{R}io de {J}aneiro 2018. {V}ol. {III}. {I}nvited
  lectures}, pages 1801--1823. World Sci. Publ., Hackensack, NJ, 2018.

\bibitem{WZ10}
W.~Winter and J.~Zacharias.
\newblock The nuclear dimension of {$\rm{C}^*$}-algebras.
\newblock {\em Adv. Math.}, 224(2):461--498, 2010.

\bibitem{yeend2006topological}
T.~Yeend.
\newblock Topological higher-rank graphs and the \textnormal{C}$^*$-algebras of
  topological 1-graphs.
\newblock In {\em Operator theory, operator algebras, and applications}, volume
  414 of {\em Contemp. Math.}, pages 231--244. American Mathematical Society,
  Providence, RI, 2006.

\bibitem{Zhang92}
S.~Zhang.
\newblock Certain {$\rm{C}^*$}-algebras with real rank zero and their corona
  and multiplier algebras. {I}.
\newblock {\em Pacific J. Math.}, 155(1):169--197, 1992.

\end{thebibliography}
\end{document}